\documentclass[12pt]{article}
\usepackage{amssymb, amsmath,  txfonts}
\begin{document}

\def\abstractname{\bf Abstract}
\def\dfrac{\displaystyle\frac}
\let\oldsection\section
\renewcommand\section{\setcounter{equation}{0}\oldsection}
\renewcommand\thesection{\arabic{section}}
\renewcommand\theequation{\thesection.\arabic{equation}}
\newtheorem{theorem}{\indent Theorem}[section]
\newtheorem{lemma}{\indent Lemma}[section]
\newtheorem{proposition}{\indent Proposition}[section]
\newtheorem{definition}{\indent Definition}[section]
\newtheorem{remark}{\indent Remark}[section]
\newtheorem{corollary}{\indent Corollary}[section]
\newtheorem{assumption}{\indent Assumption}[section]
\def\pd#1#2{\displaystyle\frac{\partial#1}{\partial#2}}
\def\d#1{\displaystyle\frac{d#1}{dt}}

\title{\bf Chemotaxis-Stokes system with $p$-Laplacian slow diffusion: global boundedness across the full $p>2$ range}
\author{Chunhua Jin\thanks{School of Mathematical Sciences, South China Normal
University, Guangzhou 510631, China. \texttt{jinchhua@126.com}}}
\title{\LARGE\bf   Chemotaxis-Stokes system with $p$-Laplacian slow diffusion: global boundedness across the full $p>2$ range
\thanks{This work is supported by NSFC (12671238)}
\\
\author{Chunhua Jin\thanks{
Corresponding author. Email: {\tt jinchhua@126.com}}
\\
\small (\it{School of Mathematical Sciences, South China
Normal University, }
\\
\small \it{Guangzhou, 510631, China}\rm)
}}

\date{}

\maketitle

\begin{abstract}
The global boundedness problem for three-dimensional chemotaxis--fluid systems with nonlinear diffusion presents substantial analytical challenges, particularly in the near-linear diffusion regime. Indeed, the porous-medium and $p$-Laplacian diffusion operators approach the classical Laplacian as their respective exponents tend to the linear limits $m\searrow1$ and $p\searrow2$, and the extra a priori control afforded by nonlinear structures becomes increasingly difficult to exploit.  This difficulty is reflected in the porous-medium setting, where establishing global boundedness throughout the full slow-diffusion range $m>1$ required a sequence of substantial developments spanning more than a decade (2010--2024). For the $p$-Laplacian case, despite significant progress, the global boundedness problem remained open in the near-linear range $2<p\leq p^*$, with $p^*\approx2.01247$, due to the saturation of classical Moser-type iterations as $p\searrow2$.

In this paper, we establish global bounded strong solutions for a three-dimensional chemotaxis--Stokes system with $p$-Laplacian slow diffusion and signal consumption, without logistic damping or any smallness condition on the initial data.  Our approach departs from the conventional $L^q$-iteration framework and develops a weighted second-order energy method that directly exploits the coercive structure of the nonlinear diffusion operator. A weighted Rellich identity yields second-order density estimates. To control the resulting boundary contribution on general smooth domains, we first use tangential differentiation along the boundary to lower the order of the boundary quantity, and then apply a trace-embedding estimate. Combined with a high-gradient/low-gradient decomposition, these estimates allow us to control the mixed chemotactic terms and close the coupled estimates for the density, signal, and fluid. The resulting estimates hold throughout $2<p\leq3$, thereby resolving the previously open range $2<p\leq p^*$. Together with existing results for $p>p^*$, our work establishes global boundedness throughout the entire slow-diffusion range $p>2$.
\end{abstract}


\medskip

\section{Introduction}\label{sec:intro}

Chemotaxis, the directed movement of cells along chemical gradients, is frequently coupled with fluid transport. This occurs in numerous biological settings, most notably in suspensions of aerobic bacteria such as \emph{Bacillus subtilis} near fluid interfaces \cite{Tuval2005}. The fluid motion is governed by the Navier--Stokes equations, giving rise to the chemotaxis--Navier--Stokes system
\begin{equation}\label{eq:cs}
\left\{
\begin{aligned}
&n_t+u\cdot\nabla n=\Delta n-\chi\nabla\cdot(n\nabla c),\\
&c_t+u\cdot\nabla c-\Delta c = -cn, \\
&u_t+\nabla\pi+a u\cdot\nabla u=\Delta u+n\nabla\phi,\\
& \operatorname{div}u=0.
\end{aligned}\right.
\end{equation}
In many biological settings, fluid inertia is negligible relative to viscous forces, allowing the Navier--Stokes equations to be replaced by the Stokes system ($a=0$); we refer to the resulting system as the chemotaxis--Stokes system. The qualitative theory of \eqref{eq:cs} has been extensively developed since its introduction \cite{Lorz2010,LiuLorz2011,Winkler2012CPDE,Winkler2016AIHP,Winkler2017TAMS}. In three dimensions, global large-data weak solutions are by now well established \cite{Winkler2016AIHP,Winkler2017TAMS}. Nevertheless, the existence of globally bounded solutions for large initial data remains a central open problem already in the linear-diffusion case \cite{Winkler2012CPDE,Winkler2018JDE}.

The diffusive term in the cell-density equation of \eqref{eq:cs} is commonly modeled by Fick's law, which assumes a constant diffusivity and a flux proportional
to the density gradient, i.e.  $J=-D\nabla n$, where $n$ denotes the cell density and $D>0$ is a constant diffusivity; in the absence of proliferation and loss, conservation of cell number then gives the linear diffusion equation $n_t=D\Delta n$. This provides a natural first approximation when cell-cell interactions have a negligible influence on collective transport. In biological processes
such as wound healing, tissue invasion, and microbial population expansion, however, finite cell size, cell geometry, and mechanical interactions can substantially
influence cell dispersal. These considerations motivate the use of nonlinear diffusion laws to describe collective transport more accurately.

One such law is porous-medium diffusion, $n_t=D\Delta n^m$ with $m>1$. A mechanical interpretation is obtained by considering the redistribution of cells driven by crowding-induced pressure. Repulsive interactions between neighbouring cells generate a cellular pressure that increases with the population density, and spatial differences in this pressure drive movement towards less crowded regions. Under the Darcy-type constitutive assumption $v=-\mu\nabla P(n)$, where $v$ is the mean cell velocity and $\mu>0$ is a mobility coefficient, the population flux is $J=nv=-\mu n\nabla P(n)$; choosing the power-law pressure
\[
P(n)=\frac{Dm}{\mu(m-1)}n^{m-1},
\qquad m>1,
\]
therefore yields $J=-D\nabla n^m$. Finite-size effects have likewise motivated porous-medium-like diffusion in models of swimming bacteria \cite{DiFrancesco2010}, while exclusion-based models of interacting cells connect cell geometry with porous-medium-type continuum descriptions \cite{Simpson2011}.

The corresponding collective diffusivity,
\[
D_{\mathrm{eff}}(n)=Dmn^{m-1},
\]
increases with the density and vanishes at zero density. The law therefore represents an enhanced collective redistribution in crowded regions, rather than an increase in the intrinsic motility of individual cells. Its degeneracy at zero density additionally permits finite-speed propagation for compactly supported solutions of the pure diffusion equation, and hence population fronts with a distinct boundary---a feature that degenerate-diffusion exponents estimated from hole-closing (wound-healing) assays are designed to capture \cite{McCue2019}.

A different modelling possibility concerns a nonlinear dependence of the collective transport on the density gradient, rather than on the density itself. Cell migration through an extracellular matrix involves cortical contraction, deformation of the cell and of its nucleus, and mechanical interactions with the surrounding matrix \cite{Maxian2020,ZaidelBar2025}. Such processes motivate more general transport laws, although they do not by themselves determine a particular relation between the population flux and the density gradient: doubly nonlinear diffusions of $p$-Laplacian type, acting on a nonlinear function of the density, have been incorporated into chemotaxis models with prevention of overcrowding \cite{Bendahmane2009}, and degenerate $p$-Laplacian Keller--Segel models have been studied in \cite{CongLiu2016}; gradient-dependent power-law diffusion has also been analysed in biologically motivated PDE--ODE systems \cite{Efendiev2022}. To represent a regime in which weak density gradients produce little redistribution, whereas stronger gradients generate a disproportionately larger collective flux, one may postulate
\[
J=-K|\nabla n|^{p-2}\nabla n,
\qquad K>0,\quad p>2,
\]
which by conservation of cell number gives the $p$-Laplacian diffusion equation $n_t=K\nabla\cdot\big(|\nabla n|^{p-2}\nabla n\big)$. Here the effective diffusion factor $K|\nabla n|^{p-2}$ vanishes as the density gradient approaches zero, while the flux magnitude $K|\nabla n|^{p-1}$ increases superlinearly with the gradient magnitude. This power-law response is a phenomenological constitutive assumption for gradient-dependent collective dispersal; it is not implied by extracellular-matrix viscoelasticity alone.

Porous-medium and $p$-Laplacian diffusion thus describe distinct forms of nonlinear population transport: the former links collective spreading to the local density and degenerates at $n=0$, whereas the latter prescribes a nonlinear response of the collective flux to spatial density variations and degenerates at $\nabla n=0$. Neither equation is uniformly parabolic, so the analysis of \eqref{eq:cs} does not carry over directly. The corresponding chemotaxis--Stokes systems, obtained from \eqref{eq:cs} by replacing the linear diffusion in the first equation with one of these two laws, require a separate treatment. They read
\begin{equation}\label{eq:nl}
n_t+u\cdot\nabla n=\Delta n^m-\chi\nabla\cdot(n\nabla c)
\qquad\text{with } m>1,
\end{equation}
for porous-medium diffusion, and
\begin{equation}\label{eq:pl}
n_t+u\cdot\nabla n=\nabla\cdot\big(|\nabla n|^{p-2}\nabla n\big)-\chi\nabla\cdot(n\nabla c)
\qquad\text{with } p>2,
\end{equation}
for $p$-Laplacian diffusion, where the diffusivities $D$ and $K$ have been normalized to one and the equations for $c$ and $(u,\pi)$ are those of \eqref{eq:cs} with $a=0$.

The porous-medium case \eqref{eq:nl} illustrates the difficulty of establishing global boundedness in the presence of degenerate diffusion. After a sequence of partial results and improvements \cite{DiFrancesco2010,TaoWinkler2012,TaoWinkler2013,Winkler2015CVPDE,Winkler2018JDE}, the full slow-diffusion range $m>1$ was covered by two independent approaches developed in 2022--2024 \cite{Winkler2022ANS,Jin2024CMS}. This progression over more than a decade highlights the substantial analytical difficulties caused by the interaction between degenerate diffusion, chemotactic aggregation, and fluid transport.

A further source of difficulty is the near-linear diffusion regime. The porous-medium and $p$-Laplacian operators approach the classical Laplacian as $m\searrow1$ and $p\searrow2$, respectively. Thus, as the exponents approach these limiting values, the nonlinear features that can be exploited to obtain additional a priori control become increasingly delicate. The near-linear regime is particularly challenging for the $p$-Laplacian system: as we explain below, the integrability gain in the available Moser-type iteration degenerates as $p\searrow2$. This provides a concrete manifestation of the loss of effectiveness of the standard iteration method near the linear-diffusion limit.

For the $p$-Laplacian system \eqref{eq:pl}, progress has likewise been achieved through Moser-type iteration schemes. Tao and Li \cite{TaoLi2019,TaoLi2020} established boundedness for $p>2.133$ and then for $p>2.091$, and a refined iteration in \cite{Jin2021JDE} pushed the threshold down to $p>p^*$, where $p^*$ is determined by the divergence of the recursion
\[
r_{k+1}=\frac{2(p-1)}{3}r_k^2+\Big(p^2-\frac{10}{3}\Big)r_k+\Big(4p^2-\frac{41}{3}p+\frac{34}{3}\Big),\qquad r_0=0;
\]
explicitly, $p^*$ is the unique root of the polynomial
\begin{equation}\label{eq:pstar}
\Delta(p)=\frac19\big(9p^4-96p^3+346p^2-600p+441\big)=0
\end{equation}
in $\big(2,\sqrt{\frac{13}{3}}\,\big)$, numerically $p^*\approx2.01247$. For $2<p\leq p^*$, the sequence $\{r_k\}$ converges to a finite limit, and the corresponding Moser iteration saturates before reaching the $L^\infty$ level. More precisely, the integrability exponent accessible through this scheme, namely $1+7(p-2)$, tends to $1$ as $p\searrow2$. At the limiting value $p=2$, the system reduces to the linear-diffusion case, for which large-data boundedness in three dimensions remains open. Thus, the interval $2<p\leq p^*$ is precisely the range left unresolved by the iteration method in \cite{Jin2021JDE}.

The purpose of the present paper is to close this gap. We establish global boundedness for the full range $2<p\leq p^*$, and our new second-order estimates are valid on the larger interval $(2,3]$. The proof is based on a mechanism that differs from the Moser-type iteration described above. We explain the main ideas after stating the problem and its boundary and initial conditions.

\medskip

In this paper we study the chemotaxis--Stokes system with $p$-Laplacian slow diffusion:
\begin{align}\label{eq:main}
\left\{
\begin{aligned}
&n_t+u\cdot\nabla n = \nabla\cdot\big(|\nabla n|^{p-2}\nabla n\big)-\chi\nabla\cdot(n\nabla c),\qquad && (x,t)\in Q:=\Omega\times \mathbb{R}^+,\\
&c_t+u\cdot\nabla c-\Delta c = -cn, && (x,t)\in Q,\\
&u_t+\nabla\pi = \Delta u+n\nabla\phi, && (x,t)\in Q,\\
&\operatorname{div}u = 0, && (x,t)\in Q,
\end{aligned}\right.
\end{align}
subject to the following boundary and initial conditions
\begin{equation}\label{eq:bc}
\big(|\nabla n|^{p-2}\nabla n-\chi n\nabla c\big)\cdot\nu\big|_{\partial \Omega}
=\frac{\partial c}{\partial\nu}\Big|_{\partial \Omega}=0,\qquad
u|_{\partial \Omega}=0,
\end{equation}
\begin{equation}\label{eq:ic}
n(x,0)=n_0(x),\qquad c(x,0)=c_0(x),\qquad u(x,0)=u_0(x),\qquad x\in \Omega,
\end{equation}
where $\Omega\subset \mathbb{R}^3$ is a bounded domain with smooth boundary, $p>2$, $\chi>0$ is the chemotactic sensitivity, $\phi\in W^{1,\infty}(\Omega)$ is a given potential, and $\nu$ is the outward unit normal of $\partial \Omega$.

\medskip

{$\bullet$ \bf Methodology: a weighted second-order energy framework.}

The standard approach to global boundedness is to combine entropy estimates for the density, such as estimates for $\int_\Omega n\ln n\,dx$, with Moser-type iterations that successively improve the $L^q$ integrability of the solution. For the present $p$-Laplacian system, this procedure yields an $L^\infty$ bound when $p>p^*$. However, in the range $2<p\leq p^*$, the iteration saturates at the integrability level associated with the exponent $1+7(p-2)$ and does not provide the required $L^\infty$ estimate. The main obstacle is therefore to obtain higher regularity without relying on a further iteration of the low-order integrability estimates.

Our approach departs from this conventional $L^q$-bootstrap procedure. Once the entropy estimate is established, rather than continuing to seek higher $L^q$ norms of the density through successive iterations, we derive estimates directly for its higher-order spatial derivatives. More precisely, we construct a weighted second-order functional that connects the available low-order information with Hessian estimates through the intrinsic structure of the $p$-Laplacian. This allows us to bypass the saturation of the Moser iteration and obtain the second-order estimates throughout $(2,3]$. The guiding idea is related to the second-order energy approach developed for the porous-medium case in \cite{Jin2024CMS}; however, the $p$-Laplacian structure gives rise to different coercivity and boundary terms, which require new estimates in the present setting.

The argument consists of two main parts. The first is to derive second-order estimates for the density by converting the degeneracy of the diffusion into a coercive dissipation. The second is to close the chemotactic feedback: the second-order estimate for the density generates terms involving second derivatives of the signal, while the corresponding estimate for the signal produces mixed terms involving second derivatives of both $n$ and $c$. We describe these two parts in turn.

\medskip

\textbf{(i) Second-order estimates for the density (Lemmas \ref{lem:rellich}--\ref{lem:be}).}
The starting point is the construction of the weighted second-order functional
\[
\mathcal F_\varepsilon(t):=\int_\Omega
\frac{\big(|\nabla n_\varepsilon|^2+\varepsilon\big)^{\frac p2}}{\rho}\,dx,
\qquad \rho:=n_\varepsilon+\sqrt\varepsilon,
\]
which is a regularized version of $\int_\Omega\frac{|\nabla n|^p}{n}\,dx$.
To estimate $\mathcal F_\varepsilon$, we use a weighted Rellich identity (Lemma~\ref{lem:rellich}). In combination with the vector identity
\[
\operatorname{div}\big((A_\varepsilon\cdot\nabla)A_\varepsilon\big)
=A_\varepsilon\cdot\nabla(\operatorname{div}A_\varepsilon)
+\operatorname{tr}\big((DA_\varepsilon)^2\big),
\qquad
A_\varepsilon=(|\nabla n_\varepsilon|^2+\varepsilon)^{\frac{p-2}{2}}\nabla n_\varepsilon,
\]
this identity converts first-order information into control of the Hessian quantity
\[
\int_\Omega\frac{|Z_\varepsilon|^2}{\rho}\,dx,
\qquad
Z_\varepsilon=M_\varepsilon^{1/2}D^2n_\varepsilon M_\varepsilon^{1/2}.
\]
The resulting Rellich identity also produces a boundary contribution involving the second fundamental form, which requires a separate estimate.

The dissipation structure is as follows. Lemma~\ref{lem4-2} shows that the $p$-Laplacian contribution to $\frac{d}{dt}\mathcal F_\varepsilon$ controls the two nonnegative quantities
$$
\int_\Omega\frac{|Z_\varepsilon|^2}{\rho}\,dx
\quad\text{and}\quad
\int_\Omega\frac{\big(|\nabla n_\varepsilon|^2+\varepsilon\big)^p}{\rho^3}\,dx,
$$
with an explicit positive coefficient $\alpha_p$ for all $2<p\leq3$. By the spectral bounds for the $p$-Laplacian matrix $M_\varepsilon$ (Lemmas~\ref{lem:Mv} and \ref{lem:eig}), the first quantity controls the full Hessian in the weighted sense needed below.

A further difficulty is the boundary term
$$
\int_{\partial \Omega}\frac{(A_\varepsilon\cdot\nabla)A_\varepsilon\cdot\nu}{\rho}\,ds, \qquad
A_\varepsilon=(|\nabla n_\varepsilon|^2+\varepsilon)^{\frac{p-2}{2}}\nabla n_\varepsilon,
$$
whose sign is not favorable on a general smooth domain and which cannot simply be discarded. We emphasize that the absence of a convexity assumption is also achieved in the porous-medium work \cite{Jin2024CMS}. However, the boundary term in that work was controlled by invoking an auxiliary lemma from the literature, which is not applicable to the present $p$-Laplacian boundary expression. The key point in our treatment is to differentiate tangentially along $\partial\Omega$, thereby reducing the differential order of the boundary term, and yields
$$
(A_\varepsilon\cdot\nabla)A_\varepsilon\cdot\nu=-\mathrm{II}(A_\varepsilon,A_\varepsilon), \quad \text{with $\mathrm{II}(A_\varepsilon,A_\varepsilon)=\big\langle ({A_\varepsilon}\cdot\nabla)\nu,A_\varepsilon\big\rangle$}.
$$
This order reduction is essential: it brings the boundary expression within the scope of the trace-embedding estimate, which we apply to
\[
W_\varepsilon:=\big(|\nabla n_\varepsilon|^2+\varepsilon\big)^{\frac{p-1}{2}}\rho^{-1/2}.
\]
Combining this estimate with a Sobolev estimate for $\rho$, we show that the boundary contribution can be absorbed by the two bulk dissipation terms above. Thus, although both the porous-medium work and the present paper treat general smooth domains, the boundary estimate here requires a different argument adapted to the $p$-Laplacian structure.

\medskip

\textbf{(ii) Closing the chemotactic feedback (Lemmas \ref{lem4-7}, \ref{lem:key2} and \ref{lem:2c}).}
The density estimate in (i) generates the term $\int_\Omega n_\varepsilon|\Delta c_\varepsilon|^2\,dx$. To control it, we introduce the auxiliary functionals
\[
\mathcal G_\varepsilon(t):=\int_\Omega|\nabla c_\varepsilon|^4\,dx,
\qquad
\mathcal H_\varepsilon(t):=\int_\Omega|\Delta c_\varepsilon|^2\,dx.
\]
The estimate for $\mathcal G_\varepsilon$ is direct (Lemma~\ref{lem4-7}); the fluid transport term is controlled by a Sobolev inequality, which provides the required $L^4$ bound for $\nabla c_\varepsilon$. The main difficulty lies in estimating $\mathcal H_\varepsilon$ (Lemma~\ref{lem:2c}), whose differentiation produces the mixed Hessian term
\[
\int_\Omega c_\varepsilon\Delta n_\varepsilon\Delta c_\varepsilon\,dx.
\]
The density and signal estimates therefore cannot be closed independently. The key step is Lemma~\ref{lem:key2}, where we split the mixed term into high-gradient and low-gradient regions. The high-gradient part is absorbed by the favorable density dissipation
\[
\int_\Omega\frac{|\nabla n_\varepsilon|^{2(p-2)}|D^2n_\varepsilon|^2}{\rho}\,dx
\]
obtained in (i), whereas the low-gradient part is controlled by integration by parts. After this splitting, the remaining fluid couplings
$\int_\Omega\big(|\nabla u_\varepsilon|^3+|u_\varepsilon|^6\big)\,dx$
are handled by maximal regularity for the Stokes system (Lemma~\ref{lem:fluid}), which closes the estimate for $\mathcal H_\varepsilon$.

In this way, the second-order density dissipation obtained from the weighted Rellich identity also provides the key control needed for the signal estimate. This coupling closes the two estimates simultaneously and constitutes the technical core of the paper. We now state the main result.

\begin{assumption}\label{ass:H}
$n_0\in L^\infty(\Omega)$, $\nabla n_0^{1-\frac 1p}\in L^p(\Omega)$, $c_0\in W^{2,r}(\Omega)$ for any $r>1$,
$u_0\in W^{2,r}(\Omega)\cap W^{1,r}_0(\Omega)$ for any $r>1$ with $\operatorname{div}u_0=0$,
$\phi\in W^{1,\infty}(\Omega)$, and $n_0,c_0\ge0$ in $\Omega$.
\end{assumption}

Throughout the paper, $p^*$ denotes the number introduced in \eqref{eq:pstar}, namely the unique root of $\Delta(p)=0$ in $\big(2,\sqrt{\frac{13}{3}}\,\big)$.
The theorem below completes the boundedness theory for \eqref{eq:main}--\eqref{eq:ic} in the slow-diffusion class $p>2$. For $p>p^*$ global boundedness is already contained in \cite{Jin2021JDE}, so that the genuinely new range is $2<p\le p^*$. It should be emphasized, however, that the second-order a priori estimates of Section~\ref{sec:est} are established on the whole interval $(2,3]$, which is strictly larger than $(2,p^*]$: nothing in that argument requires $p\le p^*$, only $p\le3$, and we state the estimates in this natural generality, since they may be useful in related problems. The main result can now be stated.

\begin{theorem}[Main theorem]\label{thm:main}
Let $p>2$ and assume that Assumption~\ref{ass:H} holds. Then the problem
\eqref{eq:main}--\eqref{eq:ic} admits a global bounded strong solution $(n,c,u,\pi)\in \mathcal X_1\times\mathcal X_2\times\mathcal X_3\times \mathcal X_4$,
such that for every $\beta\in(\frac34,1)$ and $r>1$,
\allowdisplaybreaks
\begin{align}
\label{1-3}
&\sup_{t\in(0,+\infty)}\left(\|n(\cdot,t)\|_{L^\infty}+\|c(\cdot,t)\|_{W^{1,\infty}}+\|u(\cdot,t)\|_{W^{1,\infty}}
+\|\nabla n(\cdot,t)\|_{L^p}+\|A^\beta u(\cdot,t)\|_{L^2}
\right)\le M_1(\beta),
\\
\label{1-4}
&\sup_{t\in (0,+\infty)}\left(\|(n_t, \, \Delta_p n,\, |\nabla n|^{p-2}D^2n)\|_{L^2(Q_1(t))}+\|\nabla n\|_{L^{2p}(Q_1(t))}+
\|(c, u)\|_{W_r^{2,1}(Q_1(t))}+\|\nabla \pi\|_{L^r(Q_1(t))}
\right)\le M_2(r),
\end{align}
where $Q_1(t)=\Omega\times (t, t+1)$, while $M_1(\beta)$ and $M_2(r)$ are positive constants depending only on $\Omega$, $\chi$, $p$, $n_0$, $c_0$, $u_0$ and on $\beta$, respectively $r$. The classes $\mathcal X_i$ ($i=1, 2, 3, 4$) are those defined in {\bf Notations} of Section~\ref{sec:approx}.
\end{theorem}

\begin{remark}[two routes to the second-order estimate]\label{rem:two}
The local-in-time $L^2$ estimates for $\Delta_p n$ and $|\nabla n|^{p-2}D^2n$ in \eqref{1-4} are obtained by different arguments in different ranges of $p$. For $2<p\leq3$ they are a consequence of the weighted second-order energy method of Section~\ref{sec:est}: Proposition~\ref{prop:key} provides the weighted quantities, and the boundedness of $n$ converts them into the unweighted ones. For $p>p^*$ they are instead deduced from the uniform local-in-time $L^2$ estimate for $n_t$ established in \cite{Jin2021JDE}, combined with an unweighted Rellich identity for the $p$-Laplacian. Since no density weight is involved, the latter route is considerably simpler, and it also covers $p>3$, where the weighted method of Section~\ref{sec:est} is not available. Both derivations are carried out in the proof of Theorem~\ref{thm:main} (Steps 1 and 2, respectively); together they cover the whole range $p>2$, because $(2,3]\cup(p^*,\infty)=(2,\infty)$.
\end{remark}

\section{The approximated problems and preliminaries}\label{sec:approx}

For any given $\varepsilon\in(0,1)$, we consider the regularized problem
\begin{align}\label{eq:approx}
\left\{
\begin{aligned}
&\partial_t n_\varepsilon+u_\varepsilon\cdot\nabla n_\varepsilon = \nabla\cdot \left(\big(|\nabla n_\varepsilon|^2+\varepsilon\big)^{\frac{p-2}{2}}\nabla n_\varepsilon\right)-\chi\nabla\cdot(n_\varepsilon\nabla c_\varepsilon), && (x,t)\in Q,\\
&\partial_t c_\varepsilon+u_\varepsilon\cdot\nabla c_\varepsilon-\Delta c_\varepsilon = -c_\varepsilon n_\varepsilon, && (x,t)\in Q,\\
&\partial_t u_\varepsilon+\nabla \pi_\varepsilon = \Delta u_\varepsilon+n_\varepsilon\nabla\phi_\varepsilon, && (x,t)\in Q,\\
&\operatorname{div}u_\varepsilon = 0, && (x,t)\in Q,\\
&\frac{\partial n_\varepsilon}{\partial\nu}\Big|_{\partial \Omega}
=\frac{\partial c_\varepsilon}{\partial\nu}\Big|_{\partial \Omega}=0,\quad
u_\varepsilon|_{\partial \Omega}=0, &&\\
&n_\varepsilon(x,0)=n^\varepsilon_0(x),\ c_\varepsilon(x,0)=c^\varepsilon_0(x),\ u_\varepsilon(x,0)=u^\varepsilon_0(x),\quad x\in \Omega,
\end{aligned}\right.
\end{align}
where $\phi_\varepsilon\in C^{1+\alpha,\frac{\alpha}{2}}(\Omega\times[0,\infty))$, and the initial values
$n^\varepsilon_0$, $c^\varepsilon_0$, $u^\varepsilon_0 \in C^{2+\alpha}(\Omega)$
serve as the smooth approximations to $\phi$ and $(n_0, c_0, u_0)$, with
$n^\varepsilon_0$ and $c^\varepsilon_0$ being non-negative, non-trivial and sufficiently smooth functions.

Let $A$ denote the Stokes operator, that is $Aw:=-P\Delta w$,
$P:L^r(\Omega)\rightarrow L^r_\sigma(\Omega)$ is the Helmholtz projection  \cite{G1}.
And  $A$ generates a bounded
analytic semigroup $\{e^{-tA}\}_{t\ge 0}$ on $L_\sigma^r$, and the solution $u$ of \eqref{eq:approx}
can be expressed as
\begin{equation}\label{u-eq}
u_\varepsilon=e^{-tA}u_{\varepsilon 0}+\int_{0}^t e^{-(t-s)A}P(n_\varepsilon(s)\nabla\varphi_\varepsilon(s))ds.
\end{equation}

In what follows, we give some notations, which will be used throughout this paper.

\medskip

{\bf Notations:}

\begin{itemize}
\item For any matrix $Q$, and vectors $x, y$, \,$\displaystyle Q (x, y):=x^TQ y=\sum_{i,j}x_iQ_{ij}y_j$;
 $|Q|$ denote the  Frobenius norm, i.e. $\displaystyle |Q|^2:=\sum_{i, j}Q_{ij}^2$.

\item $\nu$: outward unit normal; $\mathrm{II}$ denotes the second fundamental form of
$\partial \Omega$ with the convention
\begin{equation}\label{eq:sff}
\mathrm{II}(X,Y):=\langle D_X\nu,\,Y\rangle,\qquad X,Y\in T_{\partial \Omega}\ \text{(extended arbitrarily along }\partial \Omega\text{)}.
\end{equation}
With this convention $\mathrm{II}\ge0$ on convex domains.

\item  $\boldsymbol{e}=\frac{\nabla n_\varepsilon}{|\nabla n_\varepsilon|}$, $A_\varepsilon=\big(|\nabla n_\varepsilon|^2+\varepsilon\big)^{\frac{p-2}{2}}\nabla n_\varepsilon$ and
\begin{equation}\label{eq:Dpe}
\Delta_p^{\varepsilon}n_\varepsilon:=\nabla\cdot A_\varepsilon.
\end{equation}

\item The (regularized) $p$-Laplace operator matrix:
\begin{equation}\label{eq:MZ}
M_\varepsilon:=\big(|\nabla n_\varepsilon|^2+\varepsilon\big)^{\frac{p-2}{2}}I+(p-2)\big(|\nabla n_\varepsilon|^2+\varepsilon\big)^{\frac{p-4}{2}}|\nabla n_\varepsilon|^2 \boldsymbol{e}\otimes \boldsymbol{e},
\qquad
Z_\varepsilon:=M_\varepsilon^{1/2}D^2n_\varepsilon\,M_\varepsilon^{1/2}.
\end{equation}

\item $\|\cdot\|_{L^p}:=\|\cdot\|_{L^p(\Omega)}$, $\tau:=\min\{1, T^\varepsilon_{\max}/2\}$.

\item The solution classes are given by
\begin{align*}
&\mathcal X_1=\{n\in L^\infty(Q); n\ge 0, \nabla n\in L^\infty(\mathbb R^+;L^p(\Omega)),
n_t, \Delta_p n, |\nabla n|^{p-2}D^2n \in L^2_{loc}([0,\infty); L^2(\Omega))\},
\\
&\mathcal X_2=\{c\in L^\infty(\mathbb R^+; W^{1,\infty}(\Omega));  c\ge 0, c_t, \nabla^2 c\in  L^r_{loc}([0,\infty);L^r(\Omega))
\  \text{for any }\  r>1\},
\\
&\mathcal X_3=\left\{u\in L^\infty(\mathbb R^+; W^{1,\infty}(\Omega));  A^{\beta}u\in L^\infty([0,\infty); L^2(\Omega)),
 u_t, \nabla^2 u\in  L^r_{loc}([0,\infty);L^r(\Omega))\right.
 \\
& \qquad \left.\text{for any }\   \beta\in(\frac34, 1), r>1\right\},
\\
&\mathcal X_4=\{\pi; \nabla\pi\in L^r_{loc}(\mathbb R^+;L^r(\Omega)), \ \text{for any }\  r>1\}.
\end{align*}

\item $C$, $C_i$, $M$, $\tilde C$, $\hat C$: generic constants independent of $\varepsilon$ and $T^\varepsilon_{\max}$,
depending at most on $\Omega$, $\chi$, $p$, $n_0$, $c_0$, $u_0$.
\end{itemize}

\begin{remark}[smoothness]
All structural lemmas below are stated for smooth positive functions $n\in C^3(\bar \Omega)$
with $\partial_\nu n=0$; applied to the classical solution of \eqref{eq:approx}
($n_\varepsilon\in C^{2+\alpha,1+\alpha/2}$, hence $A_\varepsilon\in C^{1+\alpha}$), they hold with the
distributional interpretation of $\operatorname{div}((A_\varepsilon\cdot\nabla)A_\varepsilon)$, which is
sufficient since these identities are only used under integration against smooth test
functions. One may also insert a standard mollification and let the mollifier tend to zero;
no uniform estimate is affected.
\end{remark}

Before going further, we give the local existence lemma and some uniform estimates known for all $p>2$, which have been given in \cite{Jin2021JDE}.

\begin{lemma}[local existence; {\cite[Lemma 3.1]{Jin2021JDE}}]\label{lem:local}
Assume that $p>2$. Then there exists $T^\varepsilon_{\max}\in(0,+\infty]$ such that
\eqref{eq:approx} admits a unique classical solution $(n_\varepsilon, c_\varepsilon, u_\varepsilon, \pi_{\varepsilon})$ with $ n_\varepsilon>0,\quad c_\varepsilon>0,$
\allowdisplaybreaks
\begin{align*}
(n_\varepsilon, c_\varepsilon, u_\varepsilon)\in C^{2+\alpha, 1+\alpha/2}\big(\Omega\times (0,T^\varepsilon_{\max})\big)\cap C\big(\overline\Omega\times[0,T^\varepsilon_{\max})\big),\notag
\\
 \pi_{\varepsilon}\in C^{1+\alpha, \alpha/2}\big(\Omega\times (0,T^\varepsilon_{\max})\big)\cap C\big(\overline\Omega\times[0,T^\varepsilon_{\max})\big),
\end{align*}
such that either $T^\varepsilon_{\max}=+\infty$, or
\begin{equation}\label{eq:local}
\limsup_{t\nearrow T^\varepsilon_{\max}}\Big(\|n_\varepsilon(\cdot,t)\|_{L^\infty}
+\|c_\varepsilon(\cdot,t)\|_{W^{1,\infty}}+\|A^\beta u_\varepsilon(\cdot,t)\|_{L^2}\Big)=+\infty
\quad\text{for some }\beta\in\Big(\frac34,1\Big).
\end{equation}
\end{lemma}

By Lemma 3.4 of \cite{Jin2021JDE}, combined with the Moser iteration technique
(see, for example, Lemma 3.8 there), we further obtain the following lemma.
\begin{lemma}
\label{lem:jin3.4}
Assume that $p>2$.
Let $(n_\varepsilon, c_\varepsilon, u_\varepsilon, \pi_\varepsilon)$ be the classical
solution of \eqref{eq:approx} in $(0, T^\varepsilon_{\max})$. If there exists $\beta\ge 0$ such that
$$
\sup_{t\in(0,T^{\varepsilon}_{\max})}\|n_{\varepsilon}(\cdot, t)\|_{L^{2+\beta}} \le M
$$
for some constant $M>0$. Then there exists $\tilde M>0$ independent of $\varepsilon$ and $T^{\varepsilon}_{\max}$ such that
$$
\sup_{t\in(0, T^{\varepsilon}_{\max})}\|n_{\varepsilon}(\cdot, t)\|_{L^\infty}\le \tilde M.
$$
\end{lemma}
Consequently, the blow-up criterion in Lemma \ref{lem:local} can be relaxed to the boundedness of $\|n_\varepsilon(\cdot,t)\|_{L^2}$, and we replace Lemma \ref{lem:local}
with the following lemma.
\begin{lemma}\label{lem:localex}
Assume that $p>2$. Then there exists $T^\varepsilon_{\max}\in(0,+\infty]$ such that
\eqref{eq:approx} admits a unique classical solution $(n_\varepsilon, c_\varepsilon, u_\varepsilon, \pi_{\varepsilon})$ with $ n_\varepsilon>0,\quad c_\varepsilon>0,$
\begin{align*}
(n_\varepsilon, c_\varepsilon, u_\varepsilon)\in C^{2+\alpha, 1+\alpha/2}\big(\Omega\times (0,T^\varepsilon_{\max})\big)\cap C\big(\overline\Omega\times[0,T^\varepsilon_{\max})\big),\notag
\\
 \pi_{\varepsilon}\in C^{1+\alpha, \alpha/2}\big(\Omega\times (0,T^\varepsilon_{\max})\big)\cap C\big(\overline\Omega\times[0,T^\varepsilon_{\max})\big),
\end{align*}
such that either $T^\varepsilon_{\max}=+\infty$, or
\begin{equation}\label{eq:localex}
\limsup_{t\nearrow T^\varepsilon_{\max}}\|n_\varepsilon(\cdot,t)\|_{L^2} =+\infty.
\end{equation}
\end{lemma}

We also give the following  uniform estimates known for all $p>2$.

\begin{lemma}[maximum principle and mass; {\cite[Lemma 3.2]{Jin2021JDE}}]\label{lem:maxp}
Let $(n_\varepsilon, c_\varepsilon, u_\varepsilon, \pi_\varepsilon)$ be the classical
solution of \eqref{eq:approx} in $(0, T^\varepsilon_{\max})$.
There exist $M_c, M_n$, independent of $\varepsilon$ and $T^\varepsilon_{\max}$, such that
\[
\sup_{t\in(0,T^\varepsilon_{\max})}\|c_\varepsilon(\cdot,t)\|_{L^\infty}\le M_c,
\qquad
\sup_{t\in(0,T^\varepsilon_{\max})}\int_\Omega n_\varepsilon(x,t)\,dx\le M_n.
\]
\end{lemma}

\begin{lemma}[basic energy; {\cite[Lemma 3.3]{Jin2021JDE}}]\label{lem:energy}
Assume that $p>2$.
Let $(n_\varepsilon, c_\varepsilon, u_\varepsilon, \pi_\varepsilon)$ be the classical
solution of \eqref{eq:approx} in $(0, T^\varepsilon_{\max})$. Then there exists $C$, independent of $\varepsilon$ and $T_{\max}^{\varepsilon}$, such that
\begin{align}
&\sup_{t\in(0, T^\varepsilon_{\max})}\int_\Omega\bigg\{\frac{|\nabla c_\varepsilon|^2}{c_\varepsilon}+n_\varepsilon\ln n_\varepsilon+|u_\varepsilon|^2\bigg\}\,dx \notag\\
&\qquad+\sup_{t\in(\tau,T^\varepsilon_{\max})}\int_{t-\tau}^{t}\int_\Omega
\bigg\{c_\varepsilon|D^2\ln c_\varepsilon|^2+\frac{|\nabla n_\varepsilon|^{p}}{n_\varepsilon}+\frac{n_\varepsilon}{c_\varepsilon}|\nabla c_\varepsilon|^2
+\frac{|\nabla c_\varepsilon|^4}{c_\varepsilon^3}+|\nabla u_\varepsilon|^2\bigg\}\,dx \,ds\le C.
\label{eq:energy}
\end{align}
\end{lemma}

\begin{remark}
Note that $\tau=\min\{1, \frac{ T^\varepsilon_{\max}}2\}$. If $T^\varepsilon_{\max}\ge2$ then
$\tau=1$ is a fixed length, while if $T^\varepsilon_{\max}<2$ then every sliding interval $(t-\tau,t)\subset(0,T^\varepsilon_{\max})$ may be replaced by
the whole interval $(0,T^\varepsilon_{\max})$,  and the same inequality applies. Since the interval is no longer than
 $(0,T^\varepsilon_{\max})$,   the resulting bounds are, in either case, uniform  in $\varepsilon$ and $T^\varepsilon_{\max}$.
\end{remark}

By Lemmas 3.7, 3.9 of \cite{Jin2021JDE}, we also have

\begin{lemma}\label{lem:iterLr}
Let $(n_\varepsilon, c_\varepsilon, u_\varepsilon, \pi_\varepsilon)$ be the classical
solution of \eqref{eq:approx} in $(0, T^\varepsilon_{\max})$, and assume that $p>2$.
Then
\[
\sup_{t\in(0, T^\varepsilon_{\max})}\|n_\varepsilon(\cdot,t)\|_{L^{1+7(p-2)}} \le C,
\]
with $C$ independent of $\varepsilon$ and $T^\varepsilon_{\max}$.
\end{lemma}

Next, we give some structural lemmas as follows.

\begin{lemma}\label{lem:DA}
$D A_\varepsilon=M_\varepsilon D^2n_\varepsilon$.
\end{lemma}

{\it\bfseries Proof.}
For any multi-index $j$, a direct calculation yields
\begin{align*}
\partial_j(A_\varepsilon)_i&=\partial_j\big(\big(|\nabla n_\varepsilon|^2+\varepsilon\big)^{\frac{p-2}{2}}\partial_in_\varepsilon\big)
\\
&=(p-2)\big(|\nabla n_\varepsilon|^2+\varepsilon\big)^{\frac{p-4}{2}}\sum_{k}\partial_kn_\varepsilon\,\partial_{jk}n_\varepsilon\,\partial_in_\varepsilon+\big(|\nabla n_\varepsilon|^2+\varepsilon\big)^{\frac{p-2}{2}}\partial_{ij}n_\varepsilon
\\
&=(p-2)\big(|\nabla n_\varepsilon|^2+\varepsilon\big)^{\frac{p-4}{2}}|\nabla n_\varepsilon|^2e_i\sum_{k}e_k\,\partial_{jk}n_\varepsilon+\big(|\nabla n_\varepsilon|^2+\varepsilon\big)^{\frac{p-2}{2}}\partial_{ij}n_\varepsilon.
\end{align*}
On the other hand, by the symmetry of $D^2n_\varepsilon$, it follows that
$$
\left(\boldsymbol{e}\otimes\boldsymbol{e} D^2n_\varepsilon\right)_{ij}=\left(
\boldsymbol{e}\boldsymbol{e}^\text{T} D^2n_\varepsilon
\right)_{ij} = e_i\left(\boldsymbol{e}^\text{T} D^2n_\varepsilon\right)_{j}
= e_i\sum_{k}e_k\partial_{kj} n_\varepsilon=e_i\sum_{k}e_k\partial_{jk} n_\varepsilon.
$$
Therefore,
$$
D A_\varepsilon=\big(|\nabla n_\varepsilon|^2+\varepsilon\big)^{\frac{p-2}{2}}
\left(I+(p-2)\frac{|\nabla n_\varepsilon|^2}{|\nabla n_\varepsilon|^2+\varepsilon}\boldsymbol{e}\otimes\boldsymbol{e}\right) D^2n_\varepsilon=
M_\varepsilon D^2n_\varepsilon.
$$
The proof is completed. \hfill $\Box$

 \begin{lemma}\label{lem:Mv}
The $3\times 3$ matrix $M_\varepsilon$ is symmetric positive definite,  and $\lambda_1$, $\lambda_2$ are its eigenvalues,
where $\lambda_2$ has an algebraic multiplicity of two, with
\[
\lambda_1=\big(|\nabla n_\varepsilon|^2+\varepsilon\big)^{\frac{p-2}{2}}
\left(1+(p-2)\frac{|\nabla n_\varepsilon|^2}{|\nabla n_\varepsilon|^2+\varepsilon}\right),
\qquad
\lambda_2=\big(|\nabla n_\varepsilon|^2+\varepsilon\big)^{\frac{p-2}{2}}.
\]
such that
$$
M_\varepsilon \boldsymbol{e}=\lambda_1\boldsymbol{e},
$$
and for any $\boldsymbol{v}\bot \boldsymbol{e}$,
$$
M_\varepsilon \boldsymbol{v}=\lambda_2\boldsymbol{v}.
$$
\end{lemma}

{\it\bfseries Proof.} It is clear that $M_\varepsilon$ is symmetric. Noticing that $\left(\boldsymbol{e}\otimes\boldsymbol{e}\right)\boldsymbol{e}=\boldsymbol{e}\left(\boldsymbol{e}^T\boldsymbol{e}\right)=\boldsymbol{e}$, therefore
$$
M_\varepsilon \boldsymbol{e}=\big(|\nabla n_\varepsilon|^2+\varepsilon\big)^{\frac{p-2}{2}}
\left(1+(p-2)\frac{|\nabla n_\varepsilon|^2}{|\nabla n_\varepsilon|^2+\varepsilon}\right)\boldsymbol{e}.
$$
For any $\boldsymbol{v}\bot \boldsymbol{e}$, we note that $\left(\boldsymbol{e}\otimes\boldsymbol{e}\right)\boldsymbol{v}=\boldsymbol{e}\left(\boldsymbol{e}^T\boldsymbol{v}\right)=0$, and thus
$$
M_\varepsilon \boldsymbol{v}=\big(|\nabla n_\varepsilon|^2+\varepsilon\big)^{\frac{p-2}{2}}\boldsymbol{v}.
$$
This implies that $\big(|\nabla n_\varepsilon|^2+\varepsilon\big)^{\frac{p-2}{2}}$ is a double eigenvalue of the matrix $M_\varepsilon$. \hfill $\Box$

\begin{lemma}\label{lem:trZ}
$\,\operatorname{tr}Z_\varepsilon=\Delta_p^{\varepsilon}n_\varepsilon$. Moreover
\[
D A_\varepsilon=M_\varepsilon D^2n_\varepsilon=M_\varepsilon^{1/2}Z_\varepsilon\,M_\varepsilon^{-1/2},
\qquad
\operatorname{tr}\big((D A_\varepsilon)^2\big)=\operatorname{tr}\big(Z_\varepsilon^2\big)=|Z_\varepsilon|^2 .
\]
\end{lemma}

{\it\bfseries Proof.}
By cyclicity of the trace and Lemma~\ref{lem:DA},
\[
\operatorname{tr}Z_\varepsilon=\operatorname{tr}\big(M_\varepsilon^{1/2}D^2n_\varepsilon\,M_\varepsilon^{1/2}\big)
=\operatorname{tr}\big(M_\varepsilon D^2n_\varepsilon\big)=\operatorname{tr}(D A_\varepsilon)=\partial_i(A_\varepsilon)_i=\Delta_p^{\varepsilon}n_\varepsilon .
\]
Noting that $D^2n_\varepsilon=M_\varepsilon^{-1/2}Z_\varepsilon M_\varepsilon^{-1/2}$ and inserting into Lemma~\ref{lem:DA} gives
$$
D A_\varepsilon=M_\varepsilon^{1/2}Z_\varepsilon M_\varepsilon^{-1/2}, \quad i.e.\quad D A_\varepsilon\sim Z_\varepsilon.
$$
Since $Z_\varepsilon$ is  symmetric,  it follows that
\[
\operatorname{tr}\big((D A_\varepsilon)^2\big)=\operatorname{tr}\big(Z_\varepsilon^2\big)
=\sum_{i,j}(Z_\varepsilon)_{ij}^2=|Z_\varepsilon|^2.
\]
 \hfill $\Box$

\begin{lemma}\label{lem:eig}
$|Z_\varepsilon|^2\ge \big(|\nabla n_\varepsilon|^2+\varepsilon\big)^{p-2}|D^2n_\varepsilon|^2$.
\end{lemma}

{\it\bfseries Proof.} By Lemma \ref{lem:Mv}, diagonalize $M_\varepsilon=Q^T\Lambda Q$ with $Q$ orthogonal and
$\Lambda=\operatorname{diag}(\lambda_1,\lambda_2,\lambda_2)$,
$\lambda_1\ge \lambda_2$. Setting $\widetilde H:=QD^2n_\varepsilon Q^T$, using the cyclicity of the trace, we have
\[
|Z_\varepsilon|^2=\operatorname{tr}\big(M_\varepsilon^{1/2}D^2n_\varepsilon\,M_\varepsilon\,D^2n_\varepsilon\,M_\varepsilon^{1/2}\big)
=\operatorname{tr}\big(M_\varepsilon D^2n_\varepsilon\,M_\varepsilon D^2n_\varepsilon\big)
=\operatorname{tr}\big(\Lambda\widetilde H\,\Lambda\widetilde H\big)
=\sum_{i,j}\lambda_i\lambda_j\,\widetilde H_{ij}^{\,2}
\ge \lambda_2^2 |\widetilde H|^2 ,
\]
Since the Frobenius norm is invariant under orthogonal conjugation, we have $$|\widetilde H|^2 = |D^2n_\varepsilon|^2.$$ This completes the proof of the lemma.
 \hfill $\Box$

\begin{lemma}[Bochner-type vector identity]\label{lem:bochner}
For any $C^2$ vector field $A$,
\[
\operatorname{div}\big((A\cdot\nabla)A\big)=A\cdot\nabla(\operatorname{div}A)
+\operatorname{tr}\big((DA)^2\big).
\]
\end{lemma}

{\it\bfseries Proof.}
Let $V:=(A\cdot\nabla)A$, i.e.\ $\displaystyle V_k=\sum_{j}A_j\partial_jA_k$. A direct calculation yields
\begin{align*}
\operatorname{div} V&=\sum_{k, j}(\partial_kA_j)(\partial_jA_k)+\sum_{k,j}A_j\,\partial_{kj}A_k
\\
&=\sum_{k, j}(\partial_kA_j)(\partial_jA_k)+\sum_{j}A_j\,\partial_{j k}A_k
\\
&=\operatorname{tr}\big((DA)^2\big)+A\cdot\nabla(\operatorname{div}A).
\end{align*}
\hfill $\Box$

\section{Some a priori estimates and the proof of Theorem \ref{thm:main}}\label{sec:est}

This section contains the heart of the paper: a uniform (in $\varepsilon$)
second-order energy estimate for the regularized problems \eqref{eq:approx}.
As noted before Theorem~\ref{thm:main}, only the range $2<p\le p^*$ has to be
treated here, although the estimates below are derived in the larger class
$2<p\le 3$. Throughout this section,
$(n_\varepsilon,c_\varepsilon,u_\varepsilon,\pi_\varepsilon)$ denotes the
classical solution of \eqref{eq:approx} on $(0,T^\varepsilon_{\max})$, and
$\rho:=n_\varepsilon+\sqrt{\varepsilon}$.

The key a priori estimate of this paper is the following.

\begin{proposition}[key a priori estimate]\label{prop:key}
Assume that $2<p\le 3$, let
$(n_\varepsilon, c_\varepsilon, u_\varepsilon, \pi_\varepsilon)$ be the classical
solution of \eqref{eq:approx} in $(0, T^\varepsilon_{\max})$, and let $\rho$ be as above. Set
$\tau=\min\{1, T^\varepsilon_{\max}/2\}$. Then there exists a positive constant $C$ independent of $\varepsilon$ and $T^\varepsilon_{\max}$, such that
\begin{align}
&\sup_{0<t<T_{\max}^\varepsilon} \int_{\Omega}\left(\frac{\big(|\nabla n_\varepsilon|^2+\varepsilon\big)^{\frac{p}{2}}}{\rho}+|\nabla c_\varepsilon|^{4}+|\Delta c_\varepsilon|^{2}\right)dx\notag
\\
&+ \sup_{\tau<t<T_{\max}^\varepsilon}\int_{t-\tau}^t\int_\Omega \left(\frac{\big(|\nabla n_\varepsilon|^2+\varepsilon\big)^{p-2}|D^2n_\varepsilon|^2}{\rho}+\frac{\big(|\nabla n_\varepsilon|^2+\varepsilon\big)^{p}}{\rho^3}\right)\,dxds\notag
\\
&+\sup_{\tau<t<T_{\max}^\varepsilon}\int_{t-\tau}^t\int_{\Omega}\left(|\nabla c_\varepsilon|^{2}| D^2 c_\varepsilon|^2 +n_\varepsilon|\nabla c_\varepsilon|^{4}
+ |\nabla\Delta c_\varepsilon|^{2} + n_\varepsilon|\Delta c_\varepsilon|^2\right)\,dxds \notag
\\
\label{sum}
&\le C.
\end{align}
\end{proposition}

The remaining part of this section is devoted to the proof of
Proposition~\ref{prop:key}. The route is as follows: the weighted Rellich
identity of Lemma~\ref{lem:rellich} reveals the strict dissipativity of the
$p$-Laplacian part (Lemma~\ref{lem4-2}); after summing all contributions
(Lemmas~\ref{lem4-3} and \ref{lem:be}), the chemotactic and fluid terms are
absorbed by the second-order estimates for $c_\varepsilon$
(Lemmas~\ref{lem4-7}, \ref{lem:key2} and \ref{lem:2c}) and the Stokes
estimates (Lemma~\ref{lem:fluid}), and the resulting differential inequality
is closed by a Gronwall argument on sliding time intervals.

We first prove the following weighted Rellich identity.

\begin{lemma}[weighted Rellich identity, regularized version]\label{lem:rellich}
Let $\rho=n_\varepsilon+\sqrt{\varepsilon}$. Then
\begin{align}
\int_\Omega\frac{(\Delta_p^{\varepsilon}n_\varepsilon)^2}{\rho}\,dx
= &\int_\Omega\frac{|Z_\varepsilon|^2}{\rho}\,dx
+\int_{\partial \Omega}\frac{\mathrm{II}(A_\varepsilon,A_\varepsilon)}{\rho}\,ds
+\int_\Omega\frac{\big(|\nabla n_\varepsilon|^2+\varepsilon\big)^{\frac{p}{2}}}{\rho^{2}}\,\Delta_p^{\varepsilon}n_\varepsilon \,dx \notag
\\
&-(p-1)\int_\Omega\frac{\big(|\nabla n_\varepsilon|^2+\varepsilon\big)^{p-2}}{\rho^{2}}\,D^2n_\varepsilon(\nabla n_\varepsilon,\nabla n_\varepsilon)\,dx
-\varepsilon\int_\Omega\frac{\big(|\nabla n_\varepsilon|^2+\varepsilon\big)^{\frac{p-2}{2}}}{\rho^{2}}\,\Delta_p^{\varepsilon}n_\varepsilon \,dx \notag
\\
&+\varepsilon(p-2)\int_\Omega\frac{\big(|\nabla n_\varepsilon|^2+\varepsilon\big)^{p-3}}{\rho^{2}}\,D^2n_\varepsilon(\nabla n_\varepsilon,\nabla n_\varepsilon)\,dx,
\label{eq:rellich}
\end{align}
where $\mathrm{II}(A_\varepsilon,A_\varepsilon)=\big\langle ({A_\varepsilon}\cdot\nabla)\nu,A_\varepsilon\big\rangle$.
\end{lemma}

{\it\bfseries Proof.} Noticing that $A_\varepsilon\cdot\nu=0$ on $\partial \Omega$, and
\begin{equation}\label{eq3-10}
A_\varepsilon\cdot\nabla n_\varepsilon=\big(|\nabla n_\varepsilon|^2+\varepsilon\big)^{\frac{p-2}{2}}|\nabla n_\varepsilon|^2
=\big(|\nabla n_\varepsilon|^2+\varepsilon\big)^{\frac{p}{2}}-\varepsilon\big(|\nabla n_\varepsilon|^2+\varepsilon\big)^{\frac{p-2}{2}},
\end{equation}
integration by parts gives
\begin{align}\label{eq3-11}
\begin{aligned}
&\int_\Omega\frac{(\Delta_p^{\varepsilon}n_\varepsilon)^2}{\rho}\,dx
=-\int_\Omega A_\varepsilon\cdot\nabla\Big(\frac{\Delta_p^{\varepsilon}n_\varepsilon}{\rho}\Big)\,dx\\
&=-\int_\Omega\frac{A_\varepsilon\cdot\nabla \Delta_p^{\varepsilon}n_\varepsilon}{\rho}\,dx
+\int_\Omega\frac{(A_\varepsilon\cdot\nabla n_\varepsilon)\Delta_p^{\varepsilon}n_\varepsilon}{\rho^{2}}\,dx
\\
&=-\int_\Omega\frac{A_\varepsilon\cdot\nabla \Delta_p^{\varepsilon}n_\varepsilon}{\rho}\,dx
+\int_\Omega \big(|\nabla n_\varepsilon|^2+\varepsilon\big)^{\frac{p}{2}}\frac{\Delta_p^{\varepsilon}n_\varepsilon}{\rho^{2}}\,dx
-\varepsilon\int_\Omega \big(|\nabla n_\varepsilon|^2+\varepsilon\big)^{\frac{p-2}{2}}\frac{\Delta_p^{\varepsilon}n_\varepsilon}{\rho^{2}}\,dx.
\end{aligned}\end{align}
From Lemmas~\ref{lem:bochner} and \ref{lem:trZ}, we see that
\[
A_\varepsilon\cdot\nabla \Delta_p^{\varepsilon}n_\varepsilon=\operatorname{div}\big((A_\varepsilon\cdot\nabla)A_\varepsilon\big)-|Z_\varepsilon|^2 .
\]
Set $V_\varepsilon=(A_\varepsilon\cdot\nabla)A_\varepsilon$, hence
\begin{align}\label{eq3-12}
\begin{aligned}
&-\int_\Omega\frac{A_\varepsilon\cdot\nabla \Delta_p^{\varepsilon}n_\varepsilon}{\rho}\,dx
=-\int_\Omega\frac{\operatorname{div}V_\varepsilon}{\rho}\,dx
+\int_\Omega\frac{|Z_\varepsilon|^2}{\rho}\,dx
\\
&=\int_\Omega V_\varepsilon\cdot\nabla\Big(\frac1\rho\Big)\,dx-\int_{\partial \Omega}\frac{V_\varepsilon\cdot\nu}{\rho}\,ds+\int_\Omega\frac{|Z_\varepsilon|^2}{\rho}\,dx
\\
&=-\int_\Omega\frac{V_\varepsilon\cdot\nabla n_\varepsilon}{\rho^{2}}\,dx-\int_{\partial \Omega}\frac{V_\varepsilon\cdot\nu}{\rho}\,ds+\int_\Omega\frac{|Z_\varepsilon|^2}{\rho}\,dx.
\end{aligned}\end{align}
Since $A_\varepsilon\cdot\nu\equiv 0$ on $\partial\Omega$, it implies  that $A_\varepsilon$ is tangential, differentiating $A_\varepsilon\cdot\nu=0$ along the tangential direction $A_\varepsilon$ gives
\begin{align}\label{eq3-13}
\begin{aligned}
0&=A_\varepsilon\nabla (A_\varepsilon\cdot\nu)=\sum_{i,j}(A_\varepsilon)_j\,\partial_j\big((A_\varepsilon)_i\nu_i\big)
\\
&=\sum_{i,j}(A_\varepsilon)_j\,\partial_j(A_\varepsilon)_i\,\nu_i+\sum_{i,j}(A_\varepsilon)_i(A_\varepsilon)_j\,\partial_j\nu_i
\\
&=V_\varepsilon\cdot\nu+
\big\langle ({A_\varepsilon}\cdot\nabla)\nu,A_\varepsilon\big\rangle
\\
&=V_\varepsilon\cdot\nu+\mathrm{II}(A_\varepsilon,A_\varepsilon).
\end{aligned}\end{align}

Using \eqref{eq3-10}, we get
\begin{align}\label{eq3-14}
\begin{aligned}
V_\varepsilon\cdot\nabla n_\varepsilon
&=A_\varepsilon\cdot\nabla(A_\varepsilon\cdot\nabla n_\varepsilon)-D^2n_\varepsilon(A_\varepsilon,A_\varepsilon)\\
&=\big(|\nabla n_\varepsilon|^2+\varepsilon\big)^{\frac{p-2}{2}}\nabla n_\varepsilon\cdot\Big(p \big(|\nabla n_\varepsilon|^2+\varepsilon\big)^{\frac{p-2}{2}}-\varepsilon(p-2)\big(|\nabla n_\varepsilon|^2+\varepsilon\big)^{\frac{p-4}{2}}\Big)D^2n_\varepsilon\,\nabla n_\varepsilon
\\
&\qquad-\big(|\nabla n_\varepsilon|^2+\varepsilon\big)^{p-2}D^2n_\varepsilon(\nabla n_\varepsilon,\nabla n_\varepsilon)
\\
&=(p-1)\big(|\nabla n_\varepsilon|^2+\varepsilon\big)^{p-2}D^2n_\varepsilon(\nabla n_\varepsilon,\nabla n_\varepsilon)
-\varepsilon(p-2)\big(|\nabla n_\varepsilon|^2+\varepsilon\big)^{p-3} D^2n_\varepsilon(\nabla n_\varepsilon,\nabla n_\varepsilon).
\end{aligned}\end{align}
Substituting \eqref{eq3-13} and \eqref{eq3-14} into \eqref{eq3-12} yields
\begin{align}\label{eq3-15}
\begin{aligned}
-\int_\Omega\frac{A_\varepsilon\cdot\nabla \Delta_p^{\varepsilon}n_\varepsilon}{\rho}\,dx =&-(p-1)\int_\Omega\frac{\big(|\nabla n_\varepsilon|^2+\varepsilon\big)^{p-2}D^2n_\varepsilon(\nabla n_\varepsilon,\nabla n_\varepsilon)}{\rho^{2}}\,dx+\int_{\partial \Omega}\frac{\mathrm{II}(A_\varepsilon,A_\varepsilon)}{\rho}\,ds\\
&+\int_\Omega\frac{|Z_\varepsilon|^2}{\rho}\,dx
+\varepsilon(p-2)\int_\Omega\frac{\big(|\nabla n_\varepsilon|^2+\varepsilon\big)^{p-3} D^2n_\varepsilon(\nabla n_\varepsilon,\nabla n_\varepsilon)}{\rho^{2}}\,dx.
\end{aligned}\end{align}
Together with \eqref{eq3-11}, we complete the proof of this lemma. \hfill $\Box$

\begin{lemma} \label{lem4-1}
Let $\rho=n_\varepsilon+\sqrt{\varepsilon}$. Then
\begin{align*}
\frac{d}{dt}\int_{\Omega}\frac{\big(|\nabla n_\varepsilon|^2+\varepsilon\big)^{\frac{p}{2}}}{\rho}dx=J_1+J_2+J_3+J_4,
\end{align*}
with
\begin{align*}
J_1&=-p\int_\Omega\frac{(\Delta_p^{\varepsilon}n_\varepsilon)^2}{\rho}\,dx
+(p-1)\int_\Omega\frac{\big(|\nabla n_\varepsilon|^2+\varepsilon\big)^{\frac{p}{2}}}{\rho^{2}}\Delta_p^{\varepsilon}n_\varepsilon \,dx, \\
J_2&=\chi\int_\Omega\Big(\frac{p}{\rho}\Delta_p^{\varepsilon}n_\varepsilon-\frac{(p-1)\big(|\nabla n_\varepsilon|^2+\varepsilon\big)^{\frac{p}{2}}}{\rho^{2}}\Big)
\nabla\cdot(n_\varepsilon\nabla c_\varepsilon)\,dx,  \\
J_3&=\int_\Omega\Big(\frac{p}{\rho}\Delta_p^{\varepsilon}n_\varepsilon-\frac{(p-1)\big(|\nabla n_\varepsilon|^2+\varepsilon\big)^{\frac{p}{2}}}{\rho^{2}}\Big)
\big(u_\varepsilon\cdot\nabla n_\varepsilon\big)\,dx,  \\
J_4&=-p \varepsilon\int_\Omega\frac{\big(|\nabla n_\varepsilon|^2+\varepsilon\big)^{\frac{p-2}{2}}}{\rho^{2}}\left(\Delta_p^{\varepsilon}n_\varepsilon-\chi\nabla\cdot(n_\varepsilon\nabla c_\varepsilon)-u_\varepsilon\cdot\nabla n_\varepsilon \right) \,dx.
\end{align*}
\end{lemma}

{\it\bfseries Proof.} A direct calculation, and using \eqref{eq3-10} gives
\begin{align*}
&\frac{d}{dt}\int_{\Omega}\frac{\big(|\nabla n_\varepsilon|^2+\varepsilon\big)^{\frac{p}{2}}}{\rho}dx=-p\int_{\Omega}{\rm div}\left(\frac{A_\varepsilon}{\rho}\right)\partial_t n_\varepsilon dx-\int_{\Omega}\frac{\big(|\nabla n_\varepsilon|^2+\varepsilon\big)^{\frac{p}{2}}}{\rho^2}\partial_t n_\varepsilon dx
\\
=&-p\int_{\Omega} \frac{{\rm div}A_\varepsilon}{\rho} \partial_t n_\varepsilon dx+p\int_{\Omega}\frac{A_\varepsilon\cdot\nabla n_\varepsilon}{\rho^2}\partial_t n_\varepsilon dx
-\int_{\Omega}\frac{\big(|\nabla n_\varepsilon|^2+\varepsilon\big)^{\frac{p}{2}}}{\rho^2}\partial_t n_\varepsilon dx
\\
=&-p\int_{\Omega} \frac{{\rm div}A_\varepsilon}{\rho} \partial_t n_\varepsilon dx+p\int_{\Omega}\frac{\big(|\nabla n_\varepsilon|^2+\varepsilon\big)^{\frac{p}{2}}}{\rho^2}\partial_t n_\varepsilon dx-\varepsilon p\int_{\Omega}\frac{\big(|\nabla n_\varepsilon|^2+\varepsilon\big)^{\frac{p-2}{2}}}{\rho^2}\partial_t n_\varepsilon dx
\\
&-\int_{\Omega}\frac{\big(|\nabla n_\varepsilon|^2+\varepsilon\big)^{\frac{p}{2}}}{\rho^2}\partial_t n_\varepsilon dx.
\end{align*}
Noticing that
$$
\partial_t n_\varepsilon= {\rm div}A_\varepsilon-\chi\nabla\cdot(n_\varepsilon\nabla c_\varepsilon)-u_\varepsilon\cdot\nabla n_\varepsilon,
$$
and this lemma is proved. \hfill $\Box$

$J_1$ is the (regularized) $p$-Laplacian part and carries the new dissipation; $J_2$ is the
chemotactic part; $J_3$ is the transport part; $J_4$ collects the error terms produced by
the regularization itself. In what follows, we calculate  $J_i$ ($i=1, 2, 3, 4$) respectively.

\begin{lemma} [the diffusion part $J_1$.]\label{lem4-2}
Let $\rho=n_\varepsilon+\sqrt{\varepsilon}$, and assume that $2<p<3+\sqrt{5}$. Then
\begin{align*}
J_1\le -p\int_{\partial \Omega}\frac{\mathrm{II}(A_\varepsilon,A_\varepsilon)}{\rho}\,ds-\frac{\alpha_p(p-1)^2}{p^2}\int_\Omega \frac{|Z_\varepsilon|^2}{\rho}dx
-\frac{\alpha_p}2\int_\Omega\frac{\big(|\nabla n_\varepsilon|^2+\varepsilon\big)^{p}}{\rho^3}\,dx+R^d_{\varepsilon},
\end{align*}
with
$$
\alpha_p=(2-\frac{p^3}{4(p-1)^2})=\frac{(p-2)\left(p-\big(3-\sqrt{5}\big)\right)\left(3+\sqrt{5}-p\right)}{4(p-1)^2}>0,
$$
and
\begin{align*}
R^d_{\varepsilon} \le &\left(\frac{\alpha_p}2+2p\right)\varepsilon\int_\Omega\frac{\big(|\nabla n_\varepsilon|^2+\varepsilon\big)^{p-1}}{\rho^3}+4p^2\varepsilon\int_\Omega\frac{\big(|\nabla n_\varepsilon|^2+\varepsilon\big)^{p-2}}{\rho^2}|D^2 n_\varepsilon|\,dx.
\end{align*}
\end{lemma}

{\it\bfseries Proof.} Using Lemma \ref{lem:rellich}, we derive that
\begin{align}
J_1&=-p\int_\Omega\frac{(\Delta_p^{\varepsilon}n_\varepsilon)^2}{\rho}\,dx
+(p-1)\int_\Omega\frac{\big(|\nabla n_\varepsilon|^2+\varepsilon\big)^{\frac{p}{2}}}{\rho^{2}}\Delta_p^{\varepsilon}n_\varepsilon \,dx \notag
\\
= &-p\int_\Omega\frac{|Z_\varepsilon|^2}{\rho}\,dx
-p\int_{\partial \Omega}\frac{\mathrm{II}(A_\varepsilon,A_\varepsilon)}{\rho}\,ds
-\int_\Omega\frac{\big(|\nabla n_\varepsilon|^2+\varepsilon\big)^{\frac{p}{2}}}{\rho^{2}}\,\Delta_p^{\varepsilon}n_\varepsilon \,dx \notag
\\
&+p(p-1)\int_\Omega\frac{\big(|\nabla n_\varepsilon|^2+\varepsilon\big)^{p-2}}{\rho^{2}}\,D^2n_\varepsilon(\nabla n_\varepsilon,\nabla n_\varepsilon)\,dx
+\varepsilon p\int_\Omega\frac{\big(|\nabla n_\varepsilon|^2+\varepsilon\big)^{\frac{p-2}{2}}}{\rho^{2}}\,\Delta_p^{\varepsilon}n_\varepsilon \,dx \notag
\\
&-\varepsilon p(p-2)\int_\Omega\frac{\big(|\nabla n_\varepsilon|^2+\varepsilon\big)^{p-3}}{\rho^{2}}\,D^2n_\varepsilon(\nabla n_\varepsilon,\nabla n_\varepsilon)\,dx.
\label{eq4-3}
\end{align}
We note that
$$
D^2n_\varepsilon\,(\nabla n_\varepsilon, \nabla n_\varepsilon)=|\nabla n_\varepsilon|^2D^2n_\varepsilon \,(\boldsymbol{e}, \boldsymbol{e}),
$$
\begin{align*}
&-\int_\Omega\frac{\big(|\nabla n_\varepsilon|^2+\varepsilon\big)^{\frac{p}{2}}}{\rho^{2}}\,\Delta_p^{\varepsilon}n_\varepsilon \,dx
=\int_\Omega \big(|\nabla n_\varepsilon|^2+\varepsilon\big)^{\frac{p-2}2}\nabla n_\varepsilon\cdot\nabla\Big(\frac{\big(|\nabla n_\varepsilon|^2+\varepsilon\big)^{\frac{p}{2}}}{\rho^{2}}\Big)\,dx
\\
&=p\int_\Omega\frac{\big(|\nabla n_\varepsilon|^2+\varepsilon\big)^{p-2}D^2n_\varepsilon(\nabla n_\varepsilon,\nabla n_\varepsilon)}{\rho^{2}}\,dx
-2\int_\Omega\frac{\big(|\nabla n_\varepsilon|^2+\varepsilon\big)^{p-1}|\nabla n_\varepsilon|^2}{\rho^{3}}\,dx,
\end{align*}
and
\begin{align*}
&\varepsilon p\int_\Omega\frac{\big(|\nabla n_\varepsilon|^2+\varepsilon\big)^{\frac{p-2}{2}}}{\rho^{2}}\,\Delta_p^{\varepsilon}n_\varepsilon \,dx
\\
&=-\varepsilon p(p-2)\int_\Omega\frac{\big(|\nabla n_\varepsilon|^2+\varepsilon\big)^{p-3}D^2n_\varepsilon(\nabla n_\varepsilon,\nabla n_\varepsilon)}{\rho^{2}}\,dx
+2\varepsilon p\int_\Omega\frac{\big(|\nabla n_\varepsilon|^2+\varepsilon\big)^{p-2}|\nabla n_\varepsilon|^2}{\rho^{3}}\,dx.
\end{align*}
Therefore,
\begin{align*}
J_1= &-p\int_\Omega\frac{|Z_\varepsilon|^2}{\rho}\,dx-p\int_{\partial\Omega}\frac{\mathrm{II}(A_\varepsilon,A_\varepsilon)}{\rho}\,ds
-2\int_\Omega\frac{\big(|\nabla n_\varepsilon|^2+\varepsilon\big)^{p-1}|\nabla n_\varepsilon|^2}{\rho^{3}}\,dx
\\
&+p^2\int_\Omega\frac{\big(|\nabla n_\varepsilon|^2+\varepsilon\big)^{p-2}}{\rho^{2}}\,D^2n_\varepsilon(\nabla n_\varepsilon,\nabla n_\varepsilon)\,dx
+2\varepsilon p\int_\Omega\frac{\big(|\nabla n_\varepsilon|^2+\varepsilon\big)^{p-2}|\nabla n_\varepsilon|^2}{\rho^{3}}\,dx
\\
&-2\varepsilon p(p-2)\int_\Omega\frac{\big(|\nabla n_\varepsilon|^2+\varepsilon\big)^{p-3}}{\rho^{2}}\,D^2n_\varepsilon(\nabla n_\varepsilon,\nabla n_\varepsilon)\,dx
\\
=&-p\int_{\partial \Omega}\frac{\mathrm{II}(A_\varepsilon,A_\varepsilon)}{\rho}\,ds
-p\int_\Omega\frac{1}{\rho}\Big|Z_\varepsilon-\frac{p}{2(p-1)}\frac{\big(|\nabla n_\varepsilon|^2+\varepsilon\big)^{\frac{p-1}{2}}|\nabla n_\varepsilon|}{\rho}\boldsymbol{e}\otimes \boldsymbol{e}\Big|^2\,dx
\\
&-\frac{p^2}{p-1}\int_\Omega\frac{\big(|\nabla n_\varepsilon|^2+\varepsilon\big)^{\frac{p-1}{2}}|\nabla n_\varepsilon|}{\rho^2}Z_\varepsilon\,(\boldsymbol{e}, \boldsymbol{e})\,dx
+p^2\int_\Omega\frac{\big(|\nabla n_\varepsilon|^2+\varepsilon\big)^{p-2}|\nabla n_\varepsilon|^2}{\rho^2} D^2 n_\varepsilon\,(\boldsymbol{e}, \boldsymbol{e})\,dx
\\
&-(2-\frac{p^3}{4(p-1)^2})\int_\Omega\frac{\big(|\nabla n_\varepsilon|^2+\varepsilon\big)^{p-1}|\nabla n_\varepsilon|^2}{\rho^3}\,dx
\\
&+2p\varepsilon\int_\Omega\frac{\big(|\nabla n_\varepsilon|^2+\varepsilon\big)^{p-2}|\nabla n_\varepsilon|^2}{\rho^3}\,dx- 2p(p-2)\varepsilon
\int_\Omega\frac{\big(|\nabla n_\varepsilon|^2+\varepsilon\big)^{p-3}|\nabla n_\varepsilon|^2}{\rho^2} D^2 n_\varepsilon\,(\boldsymbol{e}, \boldsymbol{e})\,dx.
\end{align*}
From Lemma \ref{lem:Mv}, we see that
$$
Z_\varepsilon\,(\boldsymbol{e}, \boldsymbol{e})=\boldsymbol{e}^T Z_\varepsilon\boldsymbol{e}=\boldsymbol{e}^T\sqrt{M_\varepsilon}D^2n_\varepsilon \sqrt{M_\varepsilon}\boldsymbol{e}
=\big(\sqrt{M_\varepsilon}\boldsymbol{e}\big)^TD^2n_\varepsilon \sqrt{M_\varepsilon}\boldsymbol{e}=\lambda_1 D^2n_\varepsilon\,(\boldsymbol{e}, \boldsymbol{e})
$$
with
$$
\lambda_1=\big(|\nabla n_\varepsilon|^2+\varepsilon\big)^{\frac{p-2}{2}}
\left(p-1-\frac{(p-2)\varepsilon}{|\nabla n_\varepsilon|^2+\varepsilon}\right).
$$
Therefore
\begin{align*}
&-\frac{p^2}{p-1}\int_\Omega\frac{\big(|\nabla n_\varepsilon|^2+\varepsilon\big)^{\frac{p-1}{2}}|\nabla n_\varepsilon|}{\rho^2}Z_\varepsilon\,(\boldsymbol{e}, \boldsymbol{e})\,dx
+p^2\int_\Omega\frac{\big(|\nabla n_\varepsilon|^2+\varepsilon\big)^{p-2}|\nabla n_\varepsilon|^2}{\rho^2} D^2 n_\varepsilon\,(\boldsymbol{e}, \boldsymbol{e})\,dx
\\
=&-p^2\varepsilon\int_\Omega\frac{\big(|\nabla n_\varepsilon|^2+\varepsilon\big)^{p-2}|\nabla n_\varepsilon|}{\rho^2(|\nabla n_\varepsilon|+\sqrt{|\nabla n_\varepsilon|^2+\varepsilon})} D^2 n_\varepsilon\,(\boldsymbol{e}, \boldsymbol{e})\,dx
\\
&+
\frac{p^2(p-2)}{p-1}\varepsilon\int_\Omega\frac{\big(|\nabla n_\varepsilon|^2+\varepsilon\big)^{p-2}|\nabla n_\varepsilon|}{\rho^2\sqrt{|\nabla n_\varepsilon|^2+\varepsilon}} D^2 n_\varepsilon\,(\boldsymbol{e}, \boldsymbol{e})\,dx.
\end{align*}
Substituting it into the identity $J_1$ gives
\begin{align}
J_1=&-p\int_{\partial \Omega}\frac{\mathrm{II}(A_\varepsilon,A_\varepsilon)}{\rho}\,ds
-p\int_\Omega\frac{1}{\rho}\Big|Z_\varepsilon-\frac{p}{2(p-1)}\frac{\big(|\nabla n_\varepsilon|^2+\varepsilon\big)^{\frac{p-1}{2}}|\nabla n_\varepsilon|}{\rho}\boldsymbol{e}\otimes \boldsymbol{e}\Big|^2\,dx \notag
\\
\label{eJ1}
&-\alpha_p\int_\Omega\frac{\big(|\nabla n_\varepsilon|^2+\varepsilon\big)^{p-1}|\nabla n_\varepsilon|^2}{\rho^3}\,dx+R_{1\varepsilon}^d
\end{align}
with
$$
\alpha_p=(2-\frac{p^3}{4(p-1)^2})
=-\frac{(p-2)(p^2-6p+4)}{4(p-1)^2}=\frac{(p-2)\left(p-\big(3-\sqrt{5}\big)\right)\left(3+\sqrt{5}-p\right)}{4(p-1)^2}>0
$$
when $2<p<3+\sqrt{5}$,  and
\begin{align*}
R_{1\varepsilon}^d=&p^2\varepsilon\int_\Omega\frac{\big(|\nabla n_\varepsilon|^2+\varepsilon\big)^{p-2}|\nabla n_\varepsilon|}{\rho^2}\left(\frac{p-2}{p-1}\frac{1}{\sqrt{|\nabla n_\varepsilon|^2+\varepsilon}}-\frac{1}{(|\nabla n_\varepsilon|+\sqrt{|\nabla n_\varepsilon|^2+\varepsilon})}\right) D^2 n_\varepsilon\,(\boldsymbol{e}, \boldsymbol{e})\,dx
\\
&+2p\varepsilon\int_\Omega\frac{\big(|\nabla n_\varepsilon|^2+\varepsilon\big)^{p-2}|\nabla n_\varepsilon|^2}{\rho^3}\,dx- 2p(p-2)\varepsilon
\int_\Omega\frac{\big(|\nabla n_\varepsilon|^2+\varepsilon\big)^{p-3}|\nabla n_\varepsilon|^2}{\rho^2} D^2 n_\varepsilon\,(\boldsymbol{e}, \boldsymbol{e})\,dx
\\
\le &4p^2\varepsilon\int_\Omega\frac{\big(|\nabla n_\varepsilon|^2+\varepsilon\big)^{p-2}}{\rho^2}|D^2 n_\varepsilon|\,dx
+2p\varepsilon\int_\Omega\frac{\big(|\nabla n_\varepsilon|^2+\varepsilon\big)^{p-2}|\nabla n_\varepsilon|^2}{\rho^3}\,dx
\\
\le &4p^2\varepsilon\int_\Omega\frac{\big(|\nabla n_\varepsilon|^2+\varepsilon\big)^{p-2}}{\rho^2}|D^2 n_\varepsilon|\,dx
+2p\varepsilon\int_\Omega\frac{\big(|\nabla n_\varepsilon|^2+\varepsilon\big)^{p-1}}{\rho^3}\,dx.
\end{align*}
On the other hand, we also note that
\begin{align*}
\frac{|Z_\varepsilon|^2}{\rho}\le\frac{2}{\rho}\left(\left|Z_\varepsilon-\frac{p}{2(p-1)}\frac{\big(|\nabla n_\varepsilon|^2+\varepsilon\big)^{\frac{p-1}{2}}|\nabla n_\varepsilon|}{\rho}\boldsymbol{e}\otimes \boldsymbol{e}\right|^2+\left|\frac{p}{2(p-1)}\frac{\big(|\nabla n_\varepsilon|^2+\varepsilon\big)^{\frac{p-1}{2}}|\nabla n_\varepsilon|}{\rho}\boldsymbol{e}\otimes \boldsymbol{e}\right|^2\right)
\\
\le \frac{2}{\rho}\left|Z_\varepsilon-\frac{p}{2(p-1)}\frac{\big(|\nabla n_\varepsilon|^2+\varepsilon\big)^{\frac{p-1}{2}}|\nabla n_\varepsilon|}{\rho}\boldsymbol{e}\otimes \boldsymbol{e}\right|^2+\frac{p^2}{2(p-1)^2}\frac{\big(|\nabla n_\varepsilon|^2+\varepsilon\big)^{p-1}|\nabla n_\varepsilon|^2}{\rho^3},
\end{align*}
that is
\begin{align*}
\frac{\alpha_p(p-1)^2}{p^2}\frac{|Z_\varepsilon|^2}{\rho}
\le \frac{2\alpha_p(p-1)^2}{p^2}\frac 1{\rho}\left|Z_\varepsilon-\frac{p}{2(p-1)}\frac{\big(|\nabla n_\varepsilon|^2+\varepsilon\big)^{\frac{p-1}{2}}|\nabla n_\varepsilon|}{\rho}\boldsymbol{e}\otimes \boldsymbol{e}\right|^2+\frac{\alpha_p}{2}\frac{\big(|\nabla n_\varepsilon|^2+\varepsilon\big)^{p-1}|\nabla n_\varepsilon|^2}{\rho^3}.
\end{align*}
It is easy to check that $\frac{2\alpha_p(p-1)^2}{p^2}<p$ when $p>2$, and combining with \eqref{eJ1}, we further
obtain that
\begin{align*}
J_1\le &-p\int_{\partial \Omega}\frac{\mathrm{II}(A_\varepsilon,A_\varepsilon)}{\rho}\,ds-\frac{\alpha_p(p-1)^2}{p^2}\int_\Omega \frac{|Z_\varepsilon|^2}{\rho}dx
-\frac{\alpha_p}2\int_\Omega\frac{\big(|\nabla n_\varepsilon|^2+\varepsilon\big)^{p-1}|\nabla n_\varepsilon|^2}{\rho^3}\,dx+R^d_{1\varepsilon}
\\
\le & -p\int_{\partial \Omega}\frac{\mathrm{II}(A_\varepsilon,A_\varepsilon)}{\rho}\,ds-\frac{\alpha_p(p-1)^2}{p^2}\int_\Omega \frac{|Z_\varepsilon|^2}{\rho}dx
-\frac{\alpha_p}2\int_\Omega\frac{\big(|\nabla n_\varepsilon|^2+\varepsilon\big)^{p}}{\rho^3}\,dx+R^d_{\varepsilon},
\end{align*}
with
$$
R^d_{\varepsilon}=\frac{\alpha_p}2\varepsilon\int_\Omega\frac{\big(|\nabla n_\varepsilon|^2+\varepsilon\big)^{p-1}}{\rho^3}+R^d_{1\varepsilon},
$$
and this lemma is proved.
\hfill $\Box$

We now estimate $\sum_{i=1}^4 J_i$, and arrive at the following lemma.
\begin{lemma} \label{lem4-3}
Let $\rho=n_\varepsilon+\sqrt{\varepsilon}$, and assume that $2<p\le 3$. Then
\begin{align*}
\sum_{i=1}^4 J_i\le &-p\int_{\partial \Omega}\frac{\mathrm{II}(A_\varepsilon,A_\varepsilon)}{\rho}\,ds-\frac{\alpha_p(p-1)^2}{2p^2}\int_\Omega \frac{\big(|\nabla n_\varepsilon|^2+\varepsilon\big)^{p-2}|D^2n_\varepsilon|^2}{\rho}dx -\frac{3\alpha_p}{16}\int_\Omega\frac{\big(|\nabla n_\varepsilon|^2+\varepsilon\big)^{p}}{\rho^3}\,dx
\\
&+\chi^2C_p\int_{\Omega}n_\varepsilon|\Delta c_\varepsilon|^2\, dx
+\tilde C_{p}\int_{\Omega}\left(\rho^{\frac{3-p}{p-1}}|\nabla c_\varepsilon|^{\frac{2p}{p-1}}+\rho^{\frac{3-p}{p-1}}|u_\varepsilon|^{\frac{2p}{p-1}}\right)\, dx
+C(p)\varepsilon^{p-\frac32},
\end{align*}
where $C_p$, $\tilde C_{p}$ and $C(p)$ are constants depending only on $p$.
\end{lemma}

{\it\bfseries Proof.} We first estimate $J_2, J_3$ and $J_4$. Recalling the definition of $J_i$ in Lemma \ref{lem4-1}, we derive that for any $\eta_1>0$,
\begin{align*}
J_2&=\chi\int_\Omega\Big(\frac{p}{\rho}\Delta_p^{\varepsilon}n_\varepsilon-\frac{(p-1)\big(|\nabla n_\varepsilon|^2+\varepsilon\big)^{\frac{p}{2}}}{\rho^{2}}\Big)
\nabla\cdot(n_\varepsilon\nabla c_\varepsilon)\,dx
 \\
\le & \eta_1\int_{\Omega}\frac{|\Delta_p^{\varepsilon}n_\varepsilon|^2}{\rho}\,dx+\eta_2\int_\Omega\frac{\big(|\nabla n_\varepsilon|^2+\varepsilon\big)^{p}}{\rho^3}\,dx
+\chi^2p^2\left(\frac{1}{\eta_1}+\frac{1}{\eta_2}\right)\int_{\Omega}\frac{
|\nabla n_\varepsilon|^2|\nabla c_\varepsilon|^2+|n_\varepsilon|^2|\Delta c_\varepsilon|^2}{\rho}\, dx.
\end{align*}
Similarly,
\begin{align*}
J_3&=\int_\Omega\Big(\frac{p}{\rho}\Delta_p^{\varepsilon}n_\varepsilon-\frac{(p-1)\big(|\nabla n_\varepsilon|^2+\varepsilon\big)^{\frac{p}{2}}}{\rho^{2}}\Big)(u_\varepsilon\nabla n_\varepsilon)\,dx,
 \\
\le & \eta_1\int_{\Omega}\frac{|\Delta_p^{\varepsilon}n_\varepsilon|^2}{\rho}\,dx+\eta_2\int_\Omega\frac{\big(|\nabla n_\varepsilon|^2+\varepsilon\big)^{p}}{\rho^3}\,dx
+\frac{p^2}2\left(\frac{1}{\eta_1}+\frac{1}{\eta_2}\right)\int_{\Omega}\frac{|\nabla n_\varepsilon|^2|u_\varepsilon|^2}{\rho}\, dx.
\end{align*}
Next, we estimate $J_4$.
\begin{align*}
J_4&=-p \varepsilon\int_\Omega\frac{\big(|\nabla n_\varepsilon|^2+\varepsilon\big)^{\frac{p-2}{2}}}{\rho^{2}}\left(\Delta_p^{\varepsilon}n_\varepsilon-\chi\nabla\cdot(n_\varepsilon\nabla c_\varepsilon)-u_\varepsilon\cdot\nabla n_\varepsilon \right) \,dx
\\
\le &  \eta_1\int_{\Omega}\frac{|\Delta_p^{\varepsilon}n_\varepsilon|^2}{\rho}\,dx
+\left(\frac{p^2}{4\eta_1}+1\right)\varepsilon^2\int_\Omega\frac{\big(|\nabla n_\varepsilon|^2+\varepsilon\big)^{p-2}}{\rho^{3}}dx
\\
&+p^2\chi^2\int_{\Omega}\frac{1}{\rho}\left(
|\nabla n_\varepsilon|^2|\nabla c_\varepsilon|^2+|n_\varepsilon|^2|\Delta c_\varepsilon|^2\right)\, dx
+\frac{p^2}{2}\int_{\Omega}\frac{1}{\rho}|\nabla n_\varepsilon|^2|u_\varepsilon|^2\, dx.
\end{align*}
Noticing that $Z_\varepsilon$ is symmetric and $\operatorname{tr}Z_\varepsilon=\Delta_p^{\varepsilon}n_\varepsilon$ (Lemma~\ref{lem:trZ}), then
$$(\Delta_p^{\varepsilon}n_\varepsilon)^2=(\operatorname{tr}Z_\varepsilon)^2\le3|Z_\varepsilon|^2.$$
Therefore,
\begin{align*}
J_2+J_3+J_4\le &  3\eta_1\int_{\Omega}\frac{|\Delta_p^{\varepsilon}n_\varepsilon|^2}{\rho}\,dx+2\eta_2\int_\Omega\frac{\big(|\nabla n_\varepsilon|^2+\varepsilon\big)^{p}}{\rho^3}\,dx
\\
&+\chi^2p^2\left(\frac{1}{\eta_1}+\frac{1}{\eta_2}+1\right)\int_{\Omega}\frac{1}{\rho}\left(
|\nabla n_\varepsilon|^2|\nabla c_\varepsilon|^2+|n_\varepsilon|^2|\Delta c_\varepsilon|^2\right)\, dx
\\
&+\frac{p^2}2\left(\frac{1}{\eta_1}+\frac{1}{\eta_2}+1\right)\int_{\Omega}\frac{|\nabla n_\varepsilon|^2|u_\varepsilon|^2}{\rho}\, dx
+\left(\frac{p^2}{4\eta_1}+1\right)\varepsilon^2\int_\Omega\frac{\big(|\nabla n_\varepsilon|^2+\varepsilon\big)^{p-2}}{\rho^{3}}dx
\\
\le &  3\eta_1\int_{\Omega}\frac{|\Delta_p^{\varepsilon}n_\varepsilon|^2}{\rho}\,dx+2\eta_2\int_\Omega\frac{\big(|\nabla n_\varepsilon|^2+\varepsilon\big)^{p}}{\rho^3}\,dx
+\eta_2\int_\Omega\frac{|\nabla n_\varepsilon|^{2p}}{\rho^3}\,dx
\\
&+C_{\eta_1, \eta_2}\int_{\Omega}\rho^{\frac{3-p}{p-1}}|\nabla c_\varepsilon|^{\frac{2p}{p-1}}+\rho^{\frac{3-p}{p-1}}|u_\varepsilon|^{\frac{2p}{p-1}}\, dx
+\chi^2p^2\left(\frac{1}{\eta_1}+\frac{1}{\eta_2}+1\right)\int_{\Omega}n_\varepsilon|\Delta c_\varepsilon|^2\, dx
\\
&+\left(\frac{p^2}{4\eta_1}+1\right)\varepsilon^2\int_\Omega\frac{\big(|\nabla n_\varepsilon|^2+\varepsilon\big)^{p-2}}{\rho^{3}}dx
\\
\le &  9\eta_1\int_{\Omega}\frac{|Z_\varepsilon|^2}{\rho}\,dx+3\eta_2\int_\Omega\frac{\big(|\nabla n_\varepsilon|^2+\varepsilon\big)^{p}}{\rho^3}\,dx
+\chi^2p^2\left(\frac{1}{\eta_1}+\frac{1}{\eta_2}+1\right)\int_{\Omega}n_\varepsilon|\Delta c_\varepsilon|^2\, dx
\\
&+C_{\eta_1, \eta_2}\int_{\Omega}\left(\rho^{\frac{3-p}{p-1}}|\nabla c_\varepsilon|^{\frac{2p}{p-1}}+\rho^{\frac{3-p}{p-1}}|u_\varepsilon|^{\frac{2p}{p-1}}\right)\, dx
+\left(\frac{p^2}{4\eta_1}+1\right)\varepsilon^2\int_\Omega\frac{\big(|\nabla n_\varepsilon|^2+\varepsilon\big)^{p-2}}{\rho^{3}}\,dx,
\end{align*}
here we have used $\rho>n_{\varepsilon}$.
Take $\eta_1=\frac{\alpha_p(p-1)^2}{36 p^2}$, $\eta_2=\frac{\alpha_p}{24}$ in the above inequality, it follows that
\begin{align}
J_2+J_3+J_4\le & \frac{\alpha_p(p-1)^2}{4 p^2}\int_{\Omega}\frac{|Z_\varepsilon|^2}{\rho}\,dx+\frac{\alpha_p}{8}\int_\Omega\frac{\big(|\nabla n_\varepsilon|^2+\varepsilon\big)^{p}}{\rho^3}\,dx
+\chi^2C_p\int_{\Omega}n_\varepsilon|\Delta c_\varepsilon|^2\, dx\notag
\\
\label{J2-J4}
&+\tilde C_{p}\int_{\Omega}\left(\rho^{\frac{3-p}{p-1}}|\nabla c_\varepsilon|^{\frac{2p}{p-1}}+\rho^{\frac{3-p}{p-1}}|u_\varepsilon|^{\frac{2p}{p-1}}\right)\, dx
+\hat C_p\varepsilon^2\int_\Omega\frac{\big(|\nabla n_\varepsilon|^2+\varepsilon\big)^{p-2}}{\rho^{3}}\,dx.
\end{align}
Using Lemma \ref{lem4-2}, we further obtain
\begin{align*}
\sum_{i=1}^4 J_i\le &-p\int_{\partial \Omega}\frac{\mathrm{II}(A_\varepsilon,A_\varepsilon)}{\rho}\,ds-\frac{3\alpha_p(p-1)^2}{4p^2}\int_\Omega \frac{|Z_\varepsilon|^2}{\rho}dx
-\frac{3\alpha_p}8\int_\Omega\frac{\big(|\nabla n_\varepsilon|^2+\varepsilon\big)^{p}}{\rho^3}\,dx
\\
&+\chi^2C_p\int_{\Omega}n_\varepsilon|\Delta c_\varepsilon|^2\, dx
+\tilde C_{p}\int_{\Omega}\left(\rho^{\frac{3-p}{p-1}}|\nabla c_\varepsilon|^{\frac{2p}{p-1}}+\rho^{\frac{3-p}{p-1}}|u_\varepsilon|^{\frac{2p}{p-1}}\right)\, dx
\\
&+\hat C_p\varepsilon^2\int_\Omega\frac{\big(|\nabla n_\varepsilon|^2+\varepsilon\big)^{p-2}}{\rho^{3}}\,dx+R^d_{\varepsilon}.
\end{align*}
Recalling lemma \ref{lem:eig},  we see that
$$
|Z_\varepsilon|^2\ge \big(|\nabla n_\varepsilon|^2+\varepsilon\big)^{p-2}|D^2n_\varepsilon|^2,
$$
hence,
\begin{align}
\sum_{i=1}^4 J_i\le &-p\int_{\partial \Omega}\frac{\mathrm{II}(A_\varepsilon,A_\varepsilon)}{\rho}\,ds-\frac{3\alpha_p(p-1)^2}{4p^2}\int_\Omega \frac{\big(|\nabla n_\varepsilon|^2+\varepsilon\big)^{p-2}|D^2n_\varepsilon|^2}{\rho}dx \notag
\\
&-\frac{3\alpha_p}8\int_\Omega\frac{\big(|\nabla n_\varepsilon|^2+\varepsilon\big)^{p}}{\rho^3}\,dx
+\chi^2C_p\int_{\Omega}n_\varepsilon|\Delta c_\varepsilon|^2\, dx
+\tilde C_{p}\int_{\Omega}\left(\rho^{\frac{3-p}{p-1}}|\nabla c_\varepsilon|^{\frac{2p}{p-1}}+\rho^{\frac{3-p}{p-1}}|u_\varepsilon|^{\frac{2p}{p-1}}\right)\, dx \notag
\\
\label{ineq-J}
&+\hat C_p\varepsilon^2\int_\Omega\frac{\big(|\nabla n_\varepsilon|^2+\varepsilon\big)^{p-2}}{\rho^{3}}\,dx+R^d_{\varepsilon}.
\end{align}
In what follows, we turn our attention to the part of the regularization error in \eqref{ineq-J}. From Lemma \ref{lem4-2}, we see that
\begin{align}
&\hat C_p\varepsilon^2\int_\Omega\frac{\big(|\nabla n_\varepsilon|^2+\varepsilon\big)^{p-2}}{\rho^{3}}\,dx+R^d_{\varepsilon}
\le \hat C_p\varepsilon^2\int_\Omega\frac{\big(|\nabla n_\varepsilon|^2+\varepsilon\big)^{p-2}}{\rho^{3}}\,dx \notag
\\
\label{err}
&+\left(\frac{\alpha_p}2+2p\right)\varepsilon\int_\Omega\frac{\big(|\nabla n_\varepsilon|^2+\varepsilon\big)^{p-1}}{\rho^3}+4p^2\varepsilon\int_\Omega\frac{\big(|\nabla n_\varepsilon|^2+\varepsilon\big)^{p-2}}{\rho^2}|D^2 n_\varepsilon|\,dx:=\sum_{k=1}^3I_{k\varepsilon}.
\end{align}
For the first  error term, by Young inequality, it is easy to obtain
\begin{align*}
I_{1\varepsilon}=
\hat C_p\varepsilon^2\int_\Omega\frac{\big(|\nabla n_\varepsilon|^2+\varepsilon\big)^{p-2}}{\rho^{3}}\,dx
\le \frac{\alpha_p}{16}\int_\Omega\frac{\big(|\nabla n_\varepsilon|^2+\varepsilon\big)^{p}}{\rho^3}\,dx+C_1(p)\int_\Omega\frac{\varepsilon^{p}}{\rho^3}\,dx.
\end{align*}
For the other two error terms, we have
\begin{align*}
I_{2\varepsilon}=\left(\frac{\alpha_p}2+2p\right)\varepsilon\int_\Omega\frac{\big(|\nabla n_\varepsilon|^2+\varepsilon\big)^{p-1}}{\rho^3}\,dx
\le \frac{\alpha_p}{16}\int_\Omega\frac{\big(|\nabla n_\varepsilon|^2+\varepsilon\big)^{p}}{\rho^3}\,dx+C_2(p)\int_\Omega\frac{\varepsilon^{p}}{\rho^3}\,dx,
\end{align*}
and
\begin{align*}
I_{3\varepsilon}=&4p^2\varepsilon\int_\Omega\frac{\big(|\nabla n_\varepsilon|^2+\varepsilon\big)^{p-2}}{\rho^2}|D^2 n_\varepsilon|\,dx
\\
\le&\frac{\alpha_p(p-1)^2}{4p^2}\int_\Omega \frac{\big(|\nabla n_\varepsilon|^2+\varepsilon\big)^{p-2}|D^2n_\varepsilon|^2}{\rho}dx+\tilde{C_3}(p)\varepsilon^2\int_\Omega \frac{\big(|\nabla n_\varepsilon|^2+\varepsilon\big)^{p-2}}{\rho^3}dx
\\
\le&\frac{\alpha_p(p-1)^2}{4p^2}\int_\Omega \frac{\big(|\nabla n_\varepsilon|^2+\varepsilon\big)^{p-2}|D^2n_\varepsilon|^2}{\rho}dx
+\frac{\alpha_p}{16}\int_\Omega\frac{\big(|\nabla n_\varepsilon|^2+\varepsilon\big)^{p}}{\rho^3}\,dx+C_3(p)\int_\Omega\frac{\varepsilon^{p}}{\rho^3}\,dx.
\end{align*}
Therefore,
\begin{align}
\sum_{k=1}^3I_{k\varepsilon}&\le\frac{\alpha_p(p-1)^2}{4p^2}\int_\Omega \frac{\big(|\nabla n_\varepsilon|^2+\varepsilon\big)^{p-2}|D^2n_\varepsilon|^2}{\rho}dx
+\frac{3\alpha_p}{16}\int_\Omega\frac{\big(|\nabla n_\varepsilon|^2+\varepsilon\big)^{p}}{\rho^3}\,dx+C_4(p)\int_\Omega\frac{\varepsilon^{p}}{\rho^3}\,dx\notag
\\
\label{err-1}
&\le\frac{\alpha_p(p-1)^2}{4p^2}\int_\Omega \frac{\big(|\nabla n_\varepsilon|^2+\varepsilon\big)^{p-2}|D^2n_\varepsilon|^2}{\rho}dx
+\frac{3\alpha_p}{16}\int_\Omega\frac{\big(|\nabla n_\varepsilon|^2+\varepsilon\big)^{p}}{\rho^3}\,dx+C(p)\varepsilon^{p-\frac32}
\end{align}
since $\rho=n_\varepsilon+\sqrt{\varepsilon}\ge \sqrt{\varepsilon}$; note that the last term tends to $0$ as $\varepsilon\to0$, because $p-\frac32>0$ whenever $p>2$. Combining \eqref{ineq-J}, \eqref{err} and \eqref{err-1}, we arrive at
\begin{align*}
\sum_{i=1}^4 J_i\le &-p\int_{\partial \Omega}\frac{\mathrm{II}(A_\varepsilon,A_\varepsilon)}{\rho}\,ds-\frac{\alpha_p(p-1)^2}{2p^2}\int_\Omega \frac{\big(|\nabla n_\varepsilon|^2+\varepsilon\big)^{p-2}|D^2n_\varepsilon|^2}{\rho}dx -\frac{3\alpha_p}{16}\int_\Omega\frac{\big(|\nabla n_\varepsilon|^2+\varepsilon\big)^{p}}{\rho^3}\,dx
\\
&+\chi^2C_p\int_{\Omega}n_\varepsilon|\Delta c_\varepsilon|^2\, dx
+\tilde C_{p}\int_{\Omega}\left(\rho^{\frac{3-p}{p-1}}|\nabla c_\varepsilon|^{\frac{2p}{p-1}}+\rho^{\frac{3-p}{p-1}}|u_\varepsilon|^{\frac{2p}{p-1}}\right)\, dx
+C(p)\varepsilon^{p-\frac32}.
\end{align*}
and this lemma is proved. \hfill $\Box$

\medskip

It is apparent that for a convex domain $\Omega$, the inequality $\mathrm{II}(A_\varepsilon,A_\varepsilon)\ge 0$ holds, so that the boundary term can be neglected directly. Nevertheless, the boundary term
 $\displaystyle -p\int_{\partial \Omega}\mathrm{II}(A_\varepsilon,A_\varepsilon)/\rho\,ds$ can be controlled even in the non-convex case, and we now establish such an estimate on a general smooth domain.

\begin{lemma}\label{lem:be}
Let $\rho=n_\varepsilon+\sqrt{\varepsilon}$, and assume that $2<p\le 3$. Then
\begin{align*}
-p\int_{\partial \Omega}\frac{\mathrm{II}(A_\varepsilon,A_\varepsilon)}{\rho}\,ds
\le \frac{\alpha_p(p-1)^2}{4p^2} \int_\Omega \frac{\big(|\nabla n_\varepsilon|^2+\varepsilon\big)^{p-2}|D^2n_\varepsilon|^2}{\rho}\,dx
+\frac{\alpha_p}{16}\int_\Omega \frac{\big(|\nabla n_\varepsilon|^2+\varepsilon\big)^{p}}{\rho^{3}}dx+C.
\end{align*}
And we further have
\begin{align*}
\sum_{i=1}^4 J_i\le &-\frac{\alpha_p(p-1)^2}{4p^2}\int_\Omega \frac{\big(|\nabla n_\varepsilon|^2+\varepsilon\big)^{p-2}|D^2n_\varepsilon|^2}{\rho}dx -\frac{\alpha_p}{8}\int_\Omega\frac{\big(|\nabla n_\varepsilon|^2+\varepsilon\big)^{p}}{\rho^3}\,dx
\\
&+\chi^2C_p\int_{\Omega}n_\varepsilon|\Delta c_\varepsilon|^2\, dx
+\tilde C_{p}\int_{\Omega}\left(\rho^{\frac{3-p}{p-1}}|\nabla c_\varepsilon|^{\frac{2p}{p-1}}+\rho^{\frac{3-p}{p-1}}|u_\varepsilon|^{\frac{2p}{p-1}}\right)\, dx
+C(p)\varepsilon^{p-\frac32}+C,
\end{align*}
where $C_p$, $\tilde C_{p}$ and $C(p)$ are constants depending only on $p$, and $C$ depends only on $p$, $\Omega$, and $\|n_0\|_{L^1}$ (recall that $\int_\Omega n_\varepsilon(\cdot,t)\,dx=\int_\Omega n_0\,dx$ for all $t$).
\end{lemma}

{\it\bfseries Proof.}  Noticing that
$$
|\mathrm{II}(A_\varepsilon,A_\varepsilon)|\le K|A_\varepsilon|^2,
$$
where $K$ depends only on the geometry of $\partial\Omega$, specially on the maximal absolute principal curvature of the boundary.
Therefore,
\begin{align*}
-p\int_{\partial \Omega}\frac{\mathrm{II}(A_\varepsilon,A_\varepsilon)}{\rho}\,ds
\le Kp \int_{\partial \Omega}\frac{|A_\varepsilon|^2}{\rho}\,ds
\le Kp\int_{\partial\Omega} \frac{\big(|\nabla n_\varepsilon|^2+\varepsilon\big)^{p-1}}{\rho}\,ds
:=Kp\int_{\partial \Omega}W_\varepsilon^2\,ds,
\end{align*}
with
\[
W_\varepsilon:=\frac{\big(|\nabla n_\varepsilon|^2+\varepsilon\big)^{\frac{p-1}2}}{\rho^{1/2}} .
\]
By the trace inequality, for any $\eta>0$
\[
Kp\int_{\partial \Omega}W_\varepsilon^2\,ds\le\eta\int_\Omega|\nabla W_\varepsilon|^2\,dx+C_{\eta}\int_\Omega W_\varepsilon^2\,dx .
\]
Direct differentiation gives
\[
|\nabla W_\varepsilon|^2\le C_p\Big(\frac{\big(|\nabla n_\varepsilon|^2+\varepsilon\big)^{p-2}|D^2n_\varepsilon|^2}{\rho}+\frac{\big(|\nabla n_\varepsilon|^2+\varepsilon\big)^{p}}{\rho^{3}}\Big).
\]
Using Young's inequality and the Sobolev interpolation inequality yields
\begin{align*}
C_\eta\int_\Omega W_\varepsilon^2\,dx
&=C_\eta\int_\Omega \rho^{\frac{2p-3}{p}}\frac{\big(|\nabla n_\varepsilon|^2+\varepsilon\big)^{p-1}}{\rho^{\frac{3p-3}p}}\, dx
\le\frac{\eta}2\int_\Omega \frac{\big(|\nabla n_\varepsilon|^2+\varepsilon\big)^{p}}{\rho^{3}}dx
+C_{1\eta}\int_\Omega \rho^{2p-3}\,dx
\\
&\le\frac{\eta}2\int_\Omega \frac{\big(|\nabla n_\varepsilon|^2+\varepsilon\big)^{p}}{\rho^{3}}dx
+\frac{\eta}2\int_\Omega\rho^{-3}|\nabla\rho|^{2p}\,dx
+C_{2\eta}
\\
&\le\eta\int_\Omega \frac{\big(|\nabla n_\varepsilon|^2+\varepsilon\big)^{p}}{\rho^{3}}dx+\tilde C_{\eta},
\end{align*}
since $\nabla\rho=\nabla n_\varepsilon$, where in the second step we used the Sobolev
interpolation inequality applied to $w:=\rho^{\frac{2p-3}{2p}}$:
\[
\int_\Omega\rho^{2p-3}\,dx=\int_\Omega w^{2p}\,dx
\le\frac12\int_\Omega|\nabla w|^{2p}\,dx+C
\le C\int_\Omega\rho^{-3}|\nabla\rho|^{2p}\,dx+C ,
\]
where $|\nabla w|^{2p}=\big(\frac{2p-3}{2p}\big)^{2p}\rho^{-3}|\nabla\rho|^{2p}$,
with $C$ depending only on $\|w\|_{L^1}^{2p}$, which is uniformly bounded since
$\|\rho\|_{L^1}\le C$ and $\frac{2p-3}{2p}<1$.
Summing up, for any $\eta>0$, we obtain
\begin{align*}
-p\int_{\partial \Omega}\frac{\mathrm{II}(A_\varepsilon,A_\varepsilon)}{\rho}\,ds
\le C_p\eta \int_\Omega \frac{\big(|\nabla n_\varepsilon|^2+\varepsilon\big)^{p-2}|D^2n_\varepsilon|^2}{\rho}\,dx
+(C_p+1)\eta\int_\Omega \frac{\big(|\nabla n_\varepsilon|^2+\varepsilon\big)^{p}}{\rho^{3}}dx+\tilde C_\eta.
\end{align*}
Take $\eta$ appropriately small, such that
$$
C_p\eta\le \frac{\alpha_p(p-1)^2}{4p^2}, \quad (C_p+1)\eta\le \frac{\alpha_p}{16}.
$$
Then the   estimate for the boundary term of this lemma is complete. Combining with Lemma \ref{lem4-3}, we complete the proof.
\hfill $\Box$

\medskip
Before going further, we introduce some inequalities \cite{Jin2021JDE}.
\begin{lemma}
\label{lem:Ineq}
Let $\Omega$ be a bounded domain of $\mathbb R^N$ with smooth boundary, and  $\varphi\in C^2(\overline \Omega)$
fulfilling $\frac{\partial\varphi}{\partial {\bf n}}=0$ on $\partial\Omega$. Then
\allowdisplaybreaks
\begin{align}
\label{ineq-1}
&\int_\Omega|\varphi|^{r+p-1+\frac{Ap}N}dx \le
C_1\|\varphi\|_{L^A}^{\frac{Ap}N}\int_\Omega |\nabla \varphi^{\frac{r+p-1}p}|^pdx+C_2\|\varphi\|_{L^1}^{r+p-1+\frac{Ap}N}, \quad \text{for any $r\ge 0$, $p\ge 2$, $A\ge 1$},
\\
\label{ineq-2}
&\int_\Omega|\nabla \varphi|^{r+2} dx\le (\sqrt N+r)^2\|\varphi\|_{L^\infty}^2\int_\Omega |\nabla \varphi|^{r-2}|D^2 \varphi|^2 dx, \quad \text{for any $r\ge 2$},
\end{align}
where, $C_1$, $C_2$ are constants independent of $\varphi$.
\end{lemma}

\begin{lemma}
\label{lem4-7} Assume that $p>2$.
Let $(n_\varepsilon, c_\varepsilon, u_\varepsilon, \pi_\varepsilon)$ be the classical
solution of \eqref{eq:approx} in $(0, T^\varepsilon_{\max})$. Then for any small constant $\delta>0$,
\begin{align}
&\frac{d}{dt}\int_\Omega |\nabla c_\varepsilon|^{4}\, dx+\int_\Omega |\nabla c_\varepsilon|^{4}\, dx
+2\int_{\Omega}|\nabla c_\varepsilon|^{2}| D^2 c_\varepsilon|^2 \,dx
+2\int_{\Omega} n_\varepsilon|\nabla c_\varepsilon|^{4}\,dx \notag
\\
\label{4-11}
&\le \delta\int_{\Omega}\frac{|\nabla n_\varepsilon|^{2p}}{\rho^3}\,dx+C\int_{\Omega}|u_\varepsilon|^6\, dx+C_\delta.
\end{align}
where $C$, $C_{\delta}$  are independent of $\varepsilon$ and $T^\varepsilon_{\max}$, and $C_{\delta}$ depends on $\delta$.
\end{lemma}

{\bfseries Proof.}
Applying $\nabla$ to the second equation of \eqref{eq:approx}, and
multiplying the resulting equation by $|\nabla c_\varepsilon|^{2}\nabla c_\varepsilon$, we obtain
\begin{align}
&\frac14\frac{d}{dt}\int_\Omega |\nabla c_\varepsilon|^{4}\, dx
+\int_{\Omega}|\nabla c_\varepsilon|^{2}| D^2 c_\varepsilon|^2 \,dx+2
\int_{\Omega}|\nabla c_\varepsilon|^{2}(\nabla|\nabla c_\varepsilon|)^2 \,dx
+\int_{\Omega} n_\varepsilon|\nabla c_\varepsilon|^{4}\,dx+\frac14\int_\Omega |\nabla c_\varepsilon|^{4}\, dx\nonumber
\\
&=-\int_{\Omega}c_\varepsilon |\nabla c_\varepsilon|^{2}\nabla c_\varepsilon\nabla n_\varepsilon \,dx
+\int_{\Omega}u_\varepsilon\nabla c_\varepsilon {\rm div}(|\nabla c_\varepsilon|^{2}\nabla c_\varepsilon)\,dx+
\frac12\int_{\partial\Omega}\frac{\partial(|\nabla c_\varepsilon|^2)}{\partial{\nu}}|\nabla c_\varepsilon|^{2}\, ds
+\frac14\int_\Omega |\nabla c_\varepsilon|^{4}\, dx\nonumber
\\
&\le \delta\int_{\Omega}\frac{|\nabla n_\varepsilon|^{2p}}{\rho^3}\,dx+
C_{1\delta}\int_{\Omega}\rho^{\frac{3}{2p-1}}|\nabla c_\varepsilon|^{\frac{6p}{2p-1}}\,dx+\frac14\int_\Omega |\nabla c_\varepsilon|^{4}\, dx
+\eta\int_{\Omega}|\nabla c_\varepsilon|^6 dx+C_\eta\int_{\Omega}|u_\varepsilon|^6 \, dx
\nonumber
\\
\label{4-12}
&+\frac14\int_{\Omega}|\nabla c_\varepsilon|^{2}| D^2 c_\varepsilon|^2 \,dx
+\frac 14\int_{\Omega}|\nabla c_\varepsilon|^{2}(\nabla|\nabla c_\varepsilon|)^2 \,dx+
K\int_{\partial\Omega}|\nabla c_\varepsilon|^4 \,ds
\end{align}
for any small $\delta, \eta>0$, and $C_\eta$ depends on $\eta$, $C_{\delta}$ depends on $\delta$, where we used
$\frac{\partial |\nabla c_\varepsilon|}{\partial\nu}\Big|_{\partial\Omega}\le K|\nabla c_\varepsilon|$ under the condition that  $\frac{\partial c_\varepsilon}{\partial\nu}\Big|_{\partial\Omega}=0$.

By the boundary trace embedding inequalities, and combining with \eqref{eq:energy},
we conclude that for any small $\eta>0$,
\begin{align*}
K\int_{\partial\Omega}|\nabla c_\varepsilon|^4 \,ds
\le \eta\|\nabla(|\nabla c_\varepsilon|^{2})\|_{L^2}^2
+C_\eta\||\nabla c_\varepsilon|^{2}\|_{L^{1}}^2
\le \eta\|\nabla(|\nabla c_\varepsilon|^{2})\|_{L^2}^2
+\hat C_\eta.
\end{align*}
Noticing that $\frac{3}{2p-1}<1$, $\frac{6p}{2p-1}<4$ since $p>2$, and $\rho=n_\varepsilon+\sqrt{\varepsilon}<n_\varepsilon+1$, therefore,
\begin{align}
&C_{1\delta}\int_{\Omega}\rho^{\frac{3}{2p-1}}|\nabla c_\varepsilon|^{\frac{6p}{2p-1}}\,dx+\frac14\int_\Omega |\nabla c_\varepsilon|^{4}\, dx\notag
\\
&\le C_{1\delta}\int_{\Omega}(\frac{1}{2C_{1\delta}}n_\varepsilon+C_{2\delta})(|\nabla c_\varepsilon|^4+1)\,dx+\frac14\int_\Omega |\nabla c_\varepsilon|^{4}\, dx\notag
\\
&\le \frac12 \int_{\Omega}n_\varepsilon |\nabla c_\varepsilon|^4\, dx +C_{1\delta}C_{2\delta}\int_{\Omega}|\nabla c_\varepsilon|^4\, dx+C_{3\delta}\notag
\\
\label{ineq4-11}
&\le \frac12 \int_{\Omega}n_\varepsilon |\nabla c_\varepsilon|^4\, dx +\eta\int_{\Omega}|\nabla c_\varepsilon|^6\, dx+C_{\delta,\eta},
\end{align}
since $\|n_{\varepsilon}\|_{L^1}$ is bounded, where we used the inequality $(a+b)^{\gamma}\le a^{\gamma}+b^{\gamma}$ for $\gamma\in (0,1)$; and $a^{q_1}\le \sigma a^{q_2}+C_{\sigma}$
 for any $\sigma>0$, when $0<q_1<q_2$.
From \eqref{ineq-2}, we also have
\begin{align}\label{ineq4-12}
2\eta\int_\Omega|\nabla  c_\varepsilon|^{6} \,dx\le 2\eta(\sqrt 3+4)^2\|c_\varepsilon\|_{L^\infty}^2\int_\Omega |\nabla  c_\varepsilon|^{2}|D^2 c_\varepsilon|^2\, dx
\le C_1\eta\int_\Omega |\nabla  c_\varepsilon|^{2}|D^2 c_\varepsilon|^2\, dx.
\end{align}
Substituting the above three inequalities into \eqref{4-12} with $\eta=\frac1{4(2+C_1)}$, it follows that
\begin{align}
&\frac14\frac{d}{dt}\int_\Omega |\nabla c_\varepsilon|^{4}\, dx+\frac14\int_\Omega |\nabla c_\varepsilon|^{4}\, dx
+\frac12\int_{\Omega}|\nabla c_\varepsilon|^{2}| D^2 c_\varepsilon|^2 \,dx
+\frac12\int_{\Omega} n_\varepsilon|\nabla c_\varepsilon|^{4}\,dx \notag
\\
\label{4-13}
&\le \delta\int_{\Omega}\frac{|\nabla n_\varepsilon|^{2p}}{\rho^3}\,dx+C\int_{\Omega}|u_\varepsilon|^6\, dx+C_\delta.
\end{align}
By the arbitrariness of $\delta$, we complete the proof. \hfill $\Box$

To obtain the estimate for $\int_{\Omega}n_\varepsilon|\Delta c_\varepsilon|^2 dx$, we need the following estimate.
\begin{lemma}
\label{lem:key2}Assume that $p>2$.
Let $(n_\varepsilon, c_\varepsilon, u_\varepsilon, \pi_\varepsilon)$ be the classical
solution of \eqref{eq:approx} in $(0, T^\varepsilon_{\max})$.  Then, for any $\eta, \delta\in(0,1)$,   there exists a  positive constant $\tilde C_{\delta,\eta}$, independent of $\varepsilon$ and $t$, such that
\[
\begin{aligned}
\left|
\int_\Omega
c_\varepsilon\Delta n_\varepsilon
\Delta c_\varepsilon\,dx
\right|
 \le\delta\int_\Omega
\rho|\Delta c_\varepsilon|^2\,dx+\delta\int_\Omega
|\nabla\Delta c_\varepsilon|^2\,dx+\frac{10}{\delta}\eta^2
\int_\Omega
\frac{|\nabla n_\varepsilon|^{2(p-2)}
      |D^2n_\varepsilon|^2}{\rho}\,dx+\tilde C_{\delta,\eta},
\end{aligned}
\]
where, $\tilde C_{\delta,\eta}$ is allowed to depend on the fixed parameters
$\eta$ and $\delta$.
\end{lemma}

{\it\bfseries Proof.} Recalling Lemma \ref{lem:maxp}, we see that for any $t>0$
 \[
\|c_\varepsilon(\cdot,t)\|_{L^\infty(\Omega)} \le M_c.
\]
For $\eta\in(0,1)$, set
\[
K_\eta:=\left(\frac{M_c}{\eta}\right)^{\frac1{p-2}}.
\]
We divide the proof into several steps.

\medskip
\noindent
\textbf{Step 1. Construction of a smooth cutoff function.}

Choose $\theta\in C^\infty([0,\infty))$ such that
\[
0\le\theta\le1,\qquad \theta(s)=0\ \text{for }s\le1,\qquad \theta(s)=1\ \text{for }s\ge2,
\]
and
\[
|\theta'(s)|<2
\qquad\text{for all } s\ge0.
\]
Define
\[
\zeta_\eta
:=
\theta\left(
\frac{|\nabla n_\varepsilon|}{K_\eta}
\right),
\qquad
\phi_\eta:=1-\zeta_\eta.
\]
Since $\theta$ is constant near the origin, the composition
$\theta(|\nabla n_\varepsilon|/K_\eta)$ is smooth even at
points where $\nabla n_\varepsilon=0$.

We introduce the  high-gradient region, the low-gradient region,
and the transition region as follows:
\[
\Omega_{H,\eta}:=\{x\in\Omega; |\nabla n_\varepsilon|\ge K_\eta\}, \quad \Omega_{L,\eta}:=\{x\in\Omega; |\nabla n_\varepsilon|\le 2K_\eta\},
\quad
\mathcal T_\eta:=\{x\in\Omega; K_\eta<|\nabla n_\varepsilon|<2K_\eta\}.
\]
It is easy to see that the cutoff functions satisfy
\[
\operatorname{supp}\zeta_\eta
\subset\Omega_{H,\eta},
\qquad
\operatorname{supp}\phi_\eta
\subset\Omega_{L,\eta},
\]
and
\[
\operatorname{supp}\nabla\phi_\eta
\subset\mathcal T_\eta.
\]
Consequently,
\[
-\int_\Omega
c_\varepsilon\Delta n_\varepsilon
\Delta c_\varepsilon\,dx=I_{1,H}+I_{1,L},
\]
where
\[
I_{1,H}
:=-\int_\Omega
\zeta_\eta c_\varepsilon
\Delta n_\varepsilon\Delta c_\varepsilon\,dx,\quad I_{1,L}:=-\int_\Omega
\phi_\eta c_\varepsilon
\Delta n_\varepsilon\Delta c_\varepsilon\,dx.
\]

\medskip
\noindent
\textbf{Step 2. Estimate of the high-gradient part.}

Using $|c_\varepsilon|\le M_c$ and Young's inequality,
we obtain, for every $\delta>0$,
\[
\begin{aligned}
|I_{1,H}|
&\le
M_c\int_{\Omega_{H,\eta}}
|\Delta n_\varepsilon|
|\Delta c_\varepsilon|\,dx\\
&\le
\frac{\delta}2\int_{\Omega_{H,\eta}}
\rho|\Delta c_\varepsilon|^2\,dx
+
\frac{M_c^2}{2\delta}
\int_{\Omega_{H,\eta}}
\frac{|\Delta n_\varepsilon|^2}{\rho}\,dx.
\end{aligned}
\]
The elementary inequality
$|\Delta n_\varepsilon|^2\le N|D^2n_\varepsilon|^2$ with $(N=3)$, and the definition of $K_\eta$ imply that
\[
|\nabla n_\varepsilon|^{2(p-2)}
\ge K_\eta^{2(p-2)}
=\frac{M_c^2}{\eta^2}, \text{on $\Omega_{H,\eta}$.}
\]
Equivalently,
\[
M_c^2
\le
\eta^2|\nabla n_\varepsilon|^{2(p-2)}.
\]
It follows that
\[
M_c^2\int_{\Omega_{H,\eta}}\frac{|\Delta n_\varepsilon|^2}{\rho}\,dx
\le 3\eta^2 \int_{\Omega_{H,\eta}}\frac{|\nabla n_\varepsilon|^{2(p-2)}|D^2n_\varepsilon|^2}{\rho}\,dx
\le 3\eta^2\int_\Omega\frac{|\nabla n_\varepsilon|^{2(p-2)}|D^2n_\varepsilon|^2}{\rho}\,dx.
\]
Therefore,
\begin{equation}
\label{es-I_{1H}}
|I_{1,H}|
\le
\frac{\delta}2\int_\Omega
\rho|\Delta c_\varepsilon|^2\,dx
+\frac{3}{2\delta}\eta^2\int_\Omega\frac{|\nabla n_\varepsilon|^{2(p-2)}|D^2n_\varepsilon|^2}{\rho}\,dx.
\end{equation}

\medskip
\noindent
\textbf{Step 3. Estimate of the low-gradient part.}

\medskip

Since $\partial_\nu n_\varepsilon=0$ on $\partial\Omega$, integration by parts over $\Omega$ yields
\begin{align}
I_{1,L}&=-\int_\Omega\phi_\eta c_\varepsilon \Delta n_\varepsilon\Delta c_\varepsilon\,dx \notag
\\
&=\int_\Omega\phi_\eta c_\varepsilon\nabla n_\varepsilon\cdot\nabla\Delta c_\varepsilon\,dx+
\int_\Omega\phi_\eta\Delta c_\varepsilon\nabla n_\varepsilon\cdot\nabla c_\varepsilon\,dx+
\int_\Omega c_\varepsilon\Delta c_\varepsilon\nabla n_\varepsilon\cdot\nabla\phi_\eta\,dx \notag
\\
\label{low-grad}
&=:J_1+J_2+J_3.
\end{align}

Noticing that $|\nabla n_\varepsilon|\le2K_\eta$ on the support of $\phi_\eta$, and  using Young's inequality, we arrive at
\begin{align}\label{es-J_1}
|J_1|\le 2M_cK_\eta \int_\Omega|\nabla\Delta c_\varepsilon|\,dx
\le \frac{\delta}2\int_\Omega
|\nabla\Delta c_\varepsilon|^2\,dx +\frac{2M_c^2K_\eta^2|\Omega|}{\delta}.
\end{align}

Again using $|\nabla n_\varepsilon|\le2K_\eta$ on the
support of $\phi_\eta$, we find that
\[
\begin{aligned}
|J_2|\le 2K_\eta\int_\Omega |\Delta c_\varepsilon| |\nabla c_\varepsilon|\,dx
\le \delta_1\int_\Omega|\Delta c_\varepsilon|^2\,dx
+\frac{1}{\delta_1}K_\eta^2\int_\Omega|\nabla c_\varepsilon|^2\,dx
\end{aligned}
\]
for any $\delta_1>0$.
The homogeneous Neumann boundary condition for
$c_\varepsilon$ implies
\[
\int_\Omega\Delta c_\varepsilon\,dx
=
\int_{\partial\Omega}\partial_\nu c_\varepsilon\,dS
=0.
\]
Consequently, the Poincar\'e inequality applied to
$\Delta c_\varepsilon$ gives
\[
\int_\Omega|\Delta c_\varepsilon|^2\,dx
\le C_{\Omega}\int_\Omega
|\nabla\Delta c_\varepsilon|^2\,dx.
\]
Recalling Lemma \ref{lem:maxp} and Lemma \ref{lem:energy}, we see that $\sup_{t}\|\nabla c_{\varepsilon}(t)\|_{L^2}$ is bounded.
Choosing $\delta_1$ sufficiently small in terms of
$\delta$ and $C_{\Omega}$, we conclude that
\begin{equation}\label{es-J_2}
|J_2| \le \frac{\delta}2\int_\Omega |\nabla\Delta c_\varepsilon|^2\,dx+C_{\delta,\eta}.
\end{equation}

Next, we estimate the cutoff-gradient term $J_3$.

Differentiating the cutoff function gives
\[
|\nabla\phi_\eta|
=\left|
\frac1{K_\eta}
\theta'\left(
\frac{|\nabla n_\varepsilon|}{K_\eta}
\right)
\nabla|\nabla n_\varepsilon|\right|\le \frac{2}{K_\eta}|D^2n_\varepsilon|.
\]
Thus,
\[
|\nabla n_\varepsilon|
|\nabla\phi_\eta|
\le 4|D^2n_\varepsilon|.
\]
since $K_\eta<|\nabla n_\varepsilon|<2K_\eta$ on $\mathcal T_\eta$.
Therefore,
\[
\begin{aligned}
|J_3|&\le M_c
\int_{\mathcal T_\eta} |\Delta c_\varepsilon||\nabla n_\varepsilon||\nabla\phi_\eta|\,dx
\le 4 M_c \int_{\mathcal T_\eta}
|D^2n_\varepsilon|
|\Delta c_\varepsilon|\,dx
\\
&\le \frac{\delta}2\int_{\mathcal T_\eta}
\rho|\Delta c_\varepsilon|^2\,dx
+
\frac{8}{\delta} M_c^2
\int_{\mathcal T_\eta}
\frac{|D^2n_\varepsilon|^2}{\rho}\,dx.
\end{aligned}
\]
On $\mathcal T_\eta$, we have
$|\nabla n_\varepsilon|>K_\eta=\left(\frac{M_c}{\eta}\right)^{\frac1{p-2}}$, and hence
\[
M_c^2
\le
\eta^2|\nabla n_\varepsilon|^{2(p-2)}.
\]
Thus,
\[
\begin{aligned}
M_c^2
\int_{\mathcal T_\eta}
\frac{|D^2n_\varepsilon|^2}{\rho}\,dx
&\le
\eta^2
\int_{\mathcal T_\eta}
\frac{|\nabla n_\varepsilon|^{2(p-2)}
      |D^2n_\varepsilon|^2}{\rho}\,dx\\
&\le\eta^2\int_\Omega\frac{|\nabla n_\varepsilon|^{2(p-2)}|D^2n_\varepsilon|^2}{\rho}\,dx.
\end{aligned}
\]
We conclude that
\begin{equation}\label{es-J_3}
|J_3|\le \frac{\delta}2\int_\Omega \rho|\Delta c_\varepsilon|^2\,dx
+\frac{8}{\delta}\eta^2\int_\Omega\frac{|\nabla n_\varepsilon|^{2(p-2)}|D^2n_\varepsilon|^2}{\rho}\,dx.
\end{equation}

Combining \eqref{low-grad}-\eqref{es-J_3}, we obtain
\begin{align}\label{es-I_{1L}}
|I_{1,L}|
\le &
\frac{\delta}2\int_\Omega
\rho|\Delta c_\varepsilon|^2\,dx +\delta\int_\Omega
|\nabla\Delta c_\varepsilon|^2\,dx +\frac{8}{\delta}\eta^2\int_\Omega\frac{|\nabla n_\varepsilon|^{2(p-2)}|D^2n_\varepsilon|^2}{\rho}\,dx
+\tilde C_{\delta,\eta}.
\end{align}

\medskip
\noindent
\textbf{Step 4. Conclusion.}

\medskip

Finally, adding \eqref{es-I_{1H}} and \eqref{es-I_{1L}} yields
$$
|I_{1,H}|+|I_{1,L}|\le\delta\int_\Omega
\rho|\Delta c_\varepsilon|^2\,dx+\delta\int_\Omega
|\nabla\Delta c_\varepsilon|^2\,dx+\frac{10}{\delta}\eta^2
\int_\Omega
\frac{|\nabla n_\varepsilon|^{2(p-2)}
      |D^2n_\varepsilon|^2}{\rho}\,dx+\tilde C_{\delta,\eta},
$$
and this lemma is proved. \hfill $\Box$

\begin{lemma}
\label{lem:2c}Assume that $p>2$.
Let $(n_\varepsilon, c_\varepsilon, u_\varepsilon, \pi_\varepsilon)$ be the classical
solution of \eqref{eq:approx} in $(0, T^\varepsilon_{\max})$. Then for any small constant $\eta, \delta>0$,
\begin{align}
&\frac12\frac{d}{dt}\int_\Omega |\Delta c_\varepsilon|^{2} dx+\frac12\int_\Omega |\Delta c_\varepsilon|^{2} dx
+ \frac38\int_{\Omega}|\nabla\Delta c_\varepsilon|^{2} dx+\frac{5}8\int_{\Omega} n_\varepsilon|\Delta c_\varepsilon|^2\,dx\notag
\\
\label{c3o}
&\le 80\eta^2\int_{\Omega} \frac{|\nabla n_{\varepsilon}|^{2(p-2)}|D^2 n_\varepsilon|^2}{\rho} dx+\delta\int_{\Omega}\frac{|\nabla n_\varepsilon|^{2p}}{\rho^3}\, dx
+ C \int_{\Omega}(|\nabla u_\varepsilon|^3+
|u_\varepsilon |^6) dx+\widehat C_{\eta,\delta},
\end{align}
where $C$, $\widehat C_{\eta,\delta}$ are positive constants independent of $\varepsilon$ and $t$, and $\widehat C_{\eta,\delta}$ depends on $\eta, \delta$.
\end{lemma}

{\bfseries Proof.}
Applying $\nabla$ to the second equation of \eqref{eq:approx},
multiplying the resulting equation by $-\nabla\Delta c_\varepsilon$ on both sides, then integrating it over $\Omega$,   we arrive at
\allowdisplaybreaks
\begin{align}
&\frac12\frac{d}{dt}\int_\Omega |\Delta c_\varepsilon|^{2} dx+\frac12\int_\Omega |\Delta c_\varepsilon|^{2} dx
+\int_{\Omega}|\nabla\Delta c_\varepsilon|^{2} dx\notag
\\
&=\frac12\int_\Omega |\Delta c_\varepsilon|^{2} dx+\int_{\Omega}\nabla(u_\varepsilon\cdot\nabla c_\varepsilon)\nabla\Delta c_\varepsilon dx
-\int_{\Omega}\Delta(n_\varepsilon c_\varepsilon)\Delta c_\varepsilon dx\notag
\\
&=\frac12\int_\Omega |\Delta c_\varepsilon|^{2} dx+ \int_{\Omega}(\nabla u_\varepsilon \nabla c_\varepsilon+u_\varepsilon  D^2c_\varepsilon )\nabla\Delta c_\varepsilon dx
-\int_{\Omega}(c_\varepsilon\Delta n_\varepsilon\Delta c_\varepsilon+n_\varepsilon|\Delta c_\varepsilon|^2
+2\nabla n_\varepsilon\nabla c_\varepsilon\Delta c_\varepsilon) dx\notag
\\
&\le - \int_{\Omega} n_\varepsilon|\Delta c_\varepsilon|^2\,dx+\frac12\int_\Omega |\Delta c_\varepsilon|^{2} dx+\frac{1}{16}\int_{\Omega}|\nabla\Delta c_\varepsilon|^{2} dx+\sigma\int_{\Omega}(|\nabla c_\varepsilon|^6+|D^2 c_\varepsilon|^3)\, dx \notag
\\
\label{eq4-21}
&+C_\sigma \int_{\Omega}(|\nabla u_\varepsilon|^3+
|u_\varepsilon |^6) dx-\int_{\Omega}(c_\varepsilon\Delta n_\varepsilon\Delta c_\varepsilon+2\nabla n_\varepsilon\nabla c_\varepsilon\Delta c_\varepsilon) \,dx.
\end{align}
By Sobolev interpolation inequality and the boundedness of $c_\varepsilon$, we have
$$
\|\Delta c_\varepsilon\|_{L^2}^{2}
\le \frac18\|\nabla\Delta c_\varepsilon\|_{L^2}^{2}+C_1\|c_\varepsilon\|_{L^\infty}^2\le \frac18\|\nabla\Delta c_\varepsilon\|_{L^2}^{2}+C_2.
$$
From Gagliardo-Nirenberg interpolation inequality and the boundedness of $c_\varepsilon$ (Lemma \ref{lem:maxp}), we obtain
$$
\| D^2 c_\varepsilon\|_{L^3}^3\le
C_3\|\nabla^3 c_\varepsilon\|_{L^2}^{2}\|c_\varepsilon\|_{L^\infty}+C_4\|c_\varepsilon\|_{L^\infty}^3
\le C_5\left(1+\|\nabla\Delta c_\varepsilon\|_{L^2}^{2}\right),
$$
and
$$
\| \nabla c_\varepsilon\|_{L^6}^6\le
C_6\|\nabla^3 c_\varepsilon\|_{L^2}^{2}\|c_\varepsilon\|_{L^\infty}+C_7\|c_\varepsilon\|_{L^\infty}^6
\le C_8(\|\nabla\Delta c_\varepsilon\|_{L^2}^{2}+1).
$$
By choosing $\sigma$ sufficiently small in the inequality \eqref{eq4-21} such that $(C_5+C_8)\sigma<\frac1{8}$, then we infer that
\begin{align}
&\frac12\frac{d}{dt}\int_\Omega |\Delta c_\varepsilon|^{2} dx+\frac12\int_\Omega |\Delta c_\varepsilon|^{2} dx
+\frac{3}4\int_{\Omega}|\nabla\Delta c_\varepsilon|^{2} dx+\int_{\Omega} n_\varepsilon|\Delta c_\varepsilon|^2\,dx \notag
\\
\label{eq4-22}
&\le C_9 \int_{\Omega}(|\nabla u_\varepsilon|^3+
|u_\varepsilon |^6) dx-\int_{\Omega}(c_\varepsilon\Delta n_\varepsilon\Delta c_\varepsilon+2\nabla n_\varepsilon\nabla c_\varepsilon\Delta c_\varepsilon) \,dx+C_{10}.
\end{align}
We denote the last two terms of the above inequality by $I_1, I_2$, that is,
$$
I_1=-\int_{\Omega}c_\varepsilon\Delta n_\varepsilon\Delta c_\varepsilon\, dx,\quad I_2= -2\int_{\Omega}\nabla n_\varepsilon\nabla c_\varepsilon\Delta c_\varepsilon\, dx.
$$

For any $\delta\in(0,1)$,
\begin{align}
&|I_2|\le 2\int_{\Omega}|\nabla n_\varepsilon\nabla c_\varepsilon\Delta c_\varepsilon|\, dx=2\int_{\Omega}\frac{|\nabla n_\varepsilon|}{\rho^{\frac{3}{2p}}}\,\rho^{\frac{1}{4}}|\nabla c_\varepsilon|\,\rho^{\frac{6-p}{4p}}|\Delta c_\varepsilon|\, dx\notag
\\
\label{I_{12}}
\le & \frac18\int_{\Omega}\rho |\nabla c_\varepsilon|^4 \, dx+\delta\int_{\Omega}\frac{|\nabla n_\varepsilon|^{2p}}{\rho^3}\, dx+C_\delta\int_{\Omega}\rho^{\frac{6-p}{3p-2}}|\Delta c_\varepsilon|^{\frac{4p}{3p-2}} \,dx.
\end{align}
Noting that $\frac{6-p}{3p-2}<1$, $\frac{4p}{3p-2}<2$ when $p>2$, therefore, for any small $\sigma>0$,
$$
C_\delta\rho^{\frac{6-p}{3p-2}}|\Delta c_\varepsilon|^{\frac{4p}{3p-2}}\le C_\delta\Big(\sigma\rho+C_{\sigma}\Big)(|\Delta c_\varepsilon|^2+1).
 $$
Take $\sigma=\frac{1}{8C_\delta}$ in the above inequality, combining with \eqref{I_{12}}, and recalling Lemma \ref{lem:maxp}, it gives
\begin{align}
|I_2|\le  & \frac14\int_{\Omega}\rho |\nabla c_\varepsilon|^4 \, dx+\delta\int_{\Omega}\frac{|\nabla n_\varepsilon|^{2p}}{\rho^3}\, dx+C_{1\delta}\int_{\Omega} |\Delta c_\varepsilon|^{2} \,dx+C_{2\delta}\notag
\\
\label{I_{13}}
\le &  \frac14\int_{\Omega}\rho |\nabla c_\varepsilon|^4 \, dx+\delta\int_{\Omega}\frac{|\nabla n_\varepsilon|^{2p}}{\rho^3}\, dx+
\frac18 \int_{\Omega} |\nabla\Delta c_\varepsilon|^{2} \,dx+C_{3\delta},
\end{align}
where, we used Sobolev imbedding inequality $C_{1\delta}\|\Delta c_\varepsilon\|_{L^2}^{2}\le \frac18\|\nabla\Delta c_\varepsilon\|_{L^2}^{2}+\tilde C_\delta\|c_\varepsilon\|_{L^\infty}^2$ and the boundedness of $\|c_\varepsilon\|_{L^\infty}$.

Take $\delta=\frac18$ in Lemma \ref{lem:key2}, we also have
$$
|I_1|\le \frac18\int_\Omega
\rho|\Delta c_\varepsilon|^2\,dx+\frac18\int_\Omega
|\nabla\Delta c_\varepsilon|^2\,dx+80\eta^2
\int_\Omega \frac{|\nabla n_\varepsilon|^{2(p-2)}|D^2n_\varepsilon|^2}{\rho}\,dx+\tilde C_{\eta}.
$$

Adding the above  two inequalities yields
\begin{align}
|I_1+I_2|&\le 80\eta^2\int_{\Omega} \frac{|\nabla n_{\varepsilon}|^{2(p-2)}|D^2 n_\varepsilon|^2}{\rho} dx+\frac38\int_{\Omega} \rho|\Delta c_\varepsilon|^2\,dx
+\frac14 \int_{\Omega}|\nabla\Delta c_\varepsilon|^2\, dx\notag
\\
\label{ES-c1}
&+\delta\int_{\Omega}\frac{|\nabla n_\varepsilon|^{2p}}{\rho^3}\, dx+ C_{\eta,\delta}^1
\end{align}
for any $\eta, \delta\in(0,1)$. Substituting \eqref{ES-c1} into \eqref{eq4-22}, and noticing that $\rho=n_\varepsilon+\sqrt{\varepsilon}<n_\varepsilon+1$, we obtain
\begin{align*}
&\frac12\frac{d}{dt}\int_\Omega |\Delta c_\varepsilon|^{2} dx+\frac12\int_\Omega |\Delta c_\varepsilon|^{2} dx
+ \frac{1}2\int_{\Omega}|\nabla\Delta c_\varepsilon|^{2} dx+\frac{5}8\int_{\Omega} n_\varepsilon|\Delta c_\varepsilon|^2\,dx
\\
\le & 80\eta^2\int_{\Omega} \frac{|\nabla n_{\varepsilon}|^{2(p-2)}|D^2 n_\varepsilon|^2}{\rho} dx+\delta\int_{\Omega}\frac{|\nabla n_\varepsilon|^{2p}}{\rho^3}\, dx
+ C_9 \int_{\Omega}(|\nabla u_\varepsilon|^3+
|u_\varepsilon |^6) dx +\frac38 \int_{\Omega} |\Delta c_\varepsilon|^2\,dx
+ C_{\eta,\delta}^2
\\
\le &80\eta^2\int_{\Omega} \frac{|\nabla n_{\varepsilon}|^{2(p-2)}|D^2 n_\varepsilon|^2}{\rho} dx+\delta\int_{\Omega}\frac{|\nabla n_\varepsilon|^{2p}}{\rho^3}\, dx
+ C_9 \int_{\Omega}(|\nabla u_\varepsilon|^3+
|u_\varepsilon |^6) dx\\
&+\frac18\int_{\Omega} |\nabla\Delta c_\varepsilon|^2\,dx+\widehat C_{\eta,\delta},
\end{align*}
which implies \eqref{c3o}, and we complete the proof of this lemma.
\hfill $\Box$

\medskip

Next, we have the following priori estimates for $u_\varepsilon$.
\begin{lemma}
\label{lem:fluid}Assume that $p>2$.
Let $(n_\varepsilon, c_\varepsilon, u_\varepsilon, \pi_\varepsilon)$ be the classical
solution of \eqref{eq:approx} in $(0, T^\varepsilon_{\max})$. Then
\begin{align}\label{4-21}
\sup_{t\in(\tau,T^\varepsilon_{\max})}\int_{t-\tau}^{t}
\|A u_\varepsilon\|_{L^{\frac{8p-9}3}}^{\frac{8p-9}3}ds\le C+C\sup_{t\in(\tau,T^\varepsilon_{\max})}\int_{t-\tau}^{t}\int_{\Omega}\frac{|\nabla n_\varepsilon|^{2p}}{\rho^3}dxds.
\end{align}
Furthermore, when $0<q<8p-9$, for any $\eta>0$, there exists a positive constant $C_\eta$ depending on $p, q, \eta$ such that
\begin{align}\label{4-22}
\sup_{t\in(\tau,T^\varepsilon_{\max})}\int_{t-\tau}^t\left(\|u_\varepsilon(s)\|_{L^{q}}^{q}+\|\nabla u_\varepsilon(s)\|_{L^3}^3\right)ds\le
\eta\sup_{t\in(\tau,T^\varepsilon_{\max})}\int_{t-\tau}^t\int_{\Omega}\frac{|\nabla n_\varepsilon|^{2p}}{\rho^3}dxds+C_{\eta}.
\end{align}
Here these constants $C$, $C_\eta$ are  independent of $\varepsilon$ and $T^\varepsilon_{\max}$.
\end{lemma}

{\bfseries Proof.}
Firstly, from Gagliardo-Nirenberg interpolation inequality and the boundedness of $\|\rho\|_{L^1}$ (Lemma \ref{lem:maxp}), we obtain
\begin{align}
\label{4-23}
\|n_\varepsilon\|_{L^{\frac{8p-9}3}}^{\frac{8p-9}3}\le \|\rho\|_{L^{\frac{8p-9}3}}^{\frac{8p-9}3}=\|\rho^{\frac{2p-3}{2p}}\|_{L^q}^q\le C_1\|\nabla \rho^{\frac{2p-3}{2p}}\|_{L^{2p}}^{2p}
\|\rho^{\frac{2p-3}{2p}}\|_{L^{\frac{2p}{2p-3}}}^{q_1}+C_2\|\rho\|_{L^1}^{\frac{8p-9}3}
\le C_3\int_{\Omega}\frac{|\nabla n_\varepsilon|^{2p}}{\rho^3}dx+C_4
\end{align}
with $q=\frac{2p(8p-9)}{3(2p-3)}$, $q_1=\frac{4p^2}{3(2p-2)}$.

Using  Lemma \ref{lem:iterLr},  we derive that
\begin{align}
\sup_{t\in(0,T^\varepsilon_{\max})}\|u_\varepsilon(\cdot, t)\|_{L^3}\le &\sup_{t\in(0,T^\varepsilon_{\max})}e^{-t}\|u_{\varepsilon 0}\|_{L^3}
+\sup_{t\in(0,T^\varepsilon_{\max})}\int_0^t \|e^{-(t-s)A}P(n_\varepsilon(s)\nabla\varphi_\varepsilon(s))\|_{L^3}ds\nonumber
\\
\le &\sup_{t\in(0,T^\varepsilon_{\max})}e^{-t}\|u_{\varepsilon 0}\|_{L^3}+
\sup_{t\in(0,T^\varepsilon_{\max})}\int_0^t e^{-\lambda(t-s)}(t-s)^{-\frac32\left(\frac{1}{1+7(p-2)}-\frac13\right)}\|n_\varepsilon(s)\nabla\varphi_\varepsilon\|_{L^{1+7(p-2)}}ds
\nonumber
\\
\le &\sup_{t\in(0,T^\varepsilon_{\max})}e^{-t}\|u_{\varepsilon 0}\|_{L^3}+\sup_{t\in(0,T^\varepsilon_{\max})}
\|n_\varepsilon\|_{L^{1+7(p-2)}}\|\nabla\varphi_\varepsilon\|_{L^\infty}
\int_0^{\infty} e^{-\lambda s}s^{-\frac32\left(\frac{1}{1+7(p-2)}-\frac13\right)}ds\nonumber
\\
\label{4-24}
\le & C.
\end{align}
By  $L^{p, q}$ theory of Stokes operator \cite{FK, HP} and \eqref{4-23}, we also have
\begin{align*}
\sup_{t\in(\tau,T^\varepsilon_{\max})}\int_{t-\tau}^{t}
\|A u_\varepsilon\|_{L^{\frac{8p-9}3}}^{\frac{8p-9}3}ds&\le C_5+C_6\sup_{t\in(\tau,T^\varepsilon_{\max})}\int_{t-\tau}^{t}
\|n_\varepsilon\|_{L^{\frac{8p-9}3}}^{\frac{8p-9}3}ds \notag
\\
&\le  C_7+C_8\sup_{t\in(\tau,T^\varepsilon_{\max})}\int_{t-\tau}^{t}\int_{\Omega}\frac{|\nabla n_\varepsilon|^{2p}}{\rho^3}dxds,
\end{align*}
and \eqref{4-21} is proved.
Recalling \eqref{4-24}, and using Gagliardo-Nirenberg interpolation inequality, we arrive at
\begin{align*}
\|\nabla u_\varepsilon\|_{L^3}^3\le C_{9}\|u_\varepsilon\|_{L^{3}}^{3(1-\beta)}\|A u_\varepsilon\|_{L^{\frac{8p-9}3}}^{3\beta}+C_{10}\|u_\varepsilon\|_{L^{3}}^3
\le C_{11}\left(1+\|A u_\varepsilon\|_{L^{\frac{8p-9}3}}^{3\beta}\right),
\end{align*}
with $3\beta=\frac{8p-9}{4(2p-3)}<\frac{8p-9}{3}$  when $p>2$.
Therefore, for any $\eta>0$,
\begin{align}
\label{gradu-1}
\int_{t-\tau}^t\|\nabla u_\varepsilon(s)\|_{L^3}^3\, ds
\le \eta\int_{t-\tau}^t\|A u_\varepsilon\|_{L^{\frac{8p-9}3}}^{\frac{8p-9}3}ds+C_{1\eta}.
\end{align}
Similarly, we also have
\begin{align*}
\|u_\varepsilon\|_{L^{8p-9}}^{8p-9}\le C_{12}\|u_\varepsilon\|_{L^{3}}^{\frac{2(8p-9)}3}\|A u_\varepsilon\|_{L^{\frac{8p-9}3}}^{{\frac{8p-9}3}}+C_{13}\|u_\varepsilon\|_{L^{3}}^{8p-9}\le
C_{14}\left(1+\|A u_\varepsilon\|_{L^{\frac{8p-9}3}}^{{\frac{8p-9}3}}\right).
\end{align*}
Therefore,  for any $\eta>0$,  when $q<8p-9$,
\begin{align}\label{u-1}
\int_{t-\tau}^t\|u_\varepsilon(s)\|_{L^{q}}^{q}ds\le  \eta\int_{t-\tau}^t\|u_\varepsilon(s)\|_{L^{8p-9}}^{8p-9}ds +C_{2\eta}
\le \eta C_{14} \int_{t-\tau}^t\|A u_\varepsilon\|_{L^{\frac{8p-9}3}}^{\frac{8p-9}3}ds+C_{3\eta}.
\end{align}
By the arbitrariness of $\eta$, combining \eqref{gradu-1}, \eqref{u-1} and \eqref{4-21}, we complete the proof of this lemma. \hfill $\Box$

\noindent{\bfseries Proof of Proposition \ref{prop:key}.}
Using Lemma \ref{lem4-1} and Lemma \ref{lem:be}, we obtain that
\begin{align}
&\frac{d}{dt}\int_{\Omega}\frac{\big(|\nabla n_\varepsilon|^2+\varepsilon\big)^{\frac{p}{2}}}{\rho}dx
+\frac{\alpha_p(p-1)^2}{4p^2}\int_\Omega \frac{\big(|\nabla n_\varepsilon|^2+\varepsilon\big)^{p-2}|D^2n_\varepsilon|^2}{\rho}dx +\frac{\alpha_p}{8}\int_\Omega\frac{\big(|\nabla n_\varepsilon|^2+\varepsilon\big)^{p}}{\rho^3}\,dx\notag
\\
\label{4-31}
\le &\chi^2C_p\int_{\Omega}n_\varepsilon|\Delta c_\varepsilon|^2\, dx
+\tilde C_{p}\int_{\Omega}\left(\rho^{\frac{3-p}{p-1}}|\nabla c_\varepsilon|^{\frac{2p}{p-1}}+\rho^{\frac{3-p}{p-1}}|u_\varepsilon|^{\frac{2p}{p-1}}\right)\, dx
+C(p)\varepsilon^{p-\frac32}+C,
\end{align}
Noticing that $0\le \frac{3-p}{p-1}<1$, $0<\frac{2p}{p-1}<4$, then similar to \eqref{ineq4-11}, \eqref{ineq4-12}, we have
\begin{equation}\label{4-32}
\tilde C_{p}\int_{\Omega}\rho^{\frac{3-p}{p-1}}|\nabla c_\varepsilon|^{\frac{2p}{p-1}}\, dx
\le\frac14\int_{\Omega}n_\varepsilon|\nabla c_\varepsilon|^4 dx+\frac14\int_{\Omega}|\nabla c_\varepsilon|^2|D^2 c_\varepsilon|^2 dx+C.
\end{equation}
Here $\rho^{\frac{3-p}{p-1}}\le n_\varepsilon+1$ and $|\nabla c_\varepsilon|^{\frac{2p}{p-1}}\le |\nabla c_\varepsilon|^4+1$ are used.

Similarly, using \eqref{4-23}, for any small $\eta>0$,
\begin{align*}
&\tilde C_{p}\int_{\Omega}\rho^{\frac{3-p}{p-1}}|u_\varepsilon|^{\frac{2p}{p-1}}\, dx
\le \int_{\Omega} (n_\varepsilon+1)(|u_\varepsilon|^4+1)\,dx
\\
\le& \tilde C_{p}\int_{\Omega}\left(n_\varepsilon|u_\varepsilon|^4+2|u_\varepsilon|^4\right) dx+C
\\
\le & \eta\int_{\Omega} |n_\varepsilon|^{\frac{8p-9}{3}} \, dx+C_\eta\int_{\Omega}|u_\varepsilon|^{\frac{8p-9}{2p-3}}\,dx+C
\\
\le & \eta \hat C\int_\Omega\frac{\big(|\nabla n_\varepsilon|^2+\varepsilon\big)^{p}}{\rho^3}\,dx +C_\eta\int_{\Omega}|u_\varepsilon|^{\frac{8p-9}{2p-3}}\,dx+C.
\end{align*}
By choosing $\eta$ appropriately small in the above inequality such that $\eta \hat C<\frac{\alpha_p}{64}$, then
\begin{equation}
\label{4-33}
\tilde C_{p}\int_{\Omega}\rho^{\frac{3-p}{p-1}}|u_\varepsilon|^{\frac{2p}{p-1}}\, dx
\le \frac{\alpha_p}{64} \int_\Omega\frac{\big(|\nabla n_\varepsilon|^2+\varepsilon\big)^{p}}{\rho^3}\,dx +C\int_{\Omega}|u_\varepsilon|^{\frac{8p-9}{2p-3}}\,dx+C.
\end{equation}
We also note that
\begin{align}\label{eq4-33}
\int_{\Omega}\frac{\big(|\nabla n_\varepsilon|^2+\varepsilon\big)^{\frac{p}{2}}}{\rho}dx\le \frac{\alpha_p}{64}\int_{\Omega}\frac{\big(|\nabla n_\varepsilon|^2+\varepsilon\big)^{p}}{|\rho|^3}dx+C\|\rho\|_{L^1}.
\end{align}
Substituting \eqref{4-32}--\eqref{eq4-33} into \eqref{4-31}, and noticing that $\|\rho\|_{L^1}$ is uniformly bounded, it gives
\begin{align}
&\frac{d}{dt}\int_{\Omega}\frac{\big(|\nabla n_\varepsilon|^2+\varepsilon\big)^{\frac{p}{2}}}{\rho}dx +\int_{\Omega}\frac{\big(|\nabla n_\varepsilon|^2+\varepsilon\big)^{\frac{p}{2}}}{\rho}dx
+\frac{\alpha_p(p-1)^2}{4p^2}\int_\Omega \frac{\big(|\nabla n_\varepsilon|^2+\varepsilon\big)^{p-2}|D^2n_\varepsilon|^2}{\rho}dx \notag
\\
&+\frac{3\alpha_p}{32}\int_\Omega\frac{\big(|\nabla n_\varepsilon|^2+\varepsilon\big)^{p}}{\rho^3}\,dx\notag
\\
\label{4-34}
\le &\chi^2C_p\int_{\Omega}n_\varepsilon|\Delta c_\varepsilon|^2\, dx+\frac14\int_{\Omega}n_\varepsilon|\nabla c_\varepsilon|^4 dx+\frac14\int_{\Omega}|\nabla c_\varepsilon|^2|D^2 c_\varepsilon|^2 \,dx
+C_1\int_{\Omega}|u_\varepsilon|^{\frac{8p-9}{2p-3}}\,dx+C_2.
\end{align}
Recalling \eqref{4-11}, and taking $\delta=\frac{\alpha_p}{64}$ in \eqref{4-11} gives
\begin{align}
&\frac{d}{dt}\int_\Omega |\nabla c_\varepsilon|^{4}\, dx+\int_\Omega |\nabla c_\varepsilon|^{4}\, dx
+2\int_{\Omega}|\nabla c_\varepsilon|^{2}| D^2 c_\varepsilon|^2 \,dx
+2\int_{\Omega} n_\varepsilon|\nabla c_\varepsilon|^{4}\,dx \notag
\\
\label{4-35}
&\le \frac{\alpha_p}{64}\int_{\Omega}\frac{|\nabla n_\varepsilon|^{2p}}{\rho^3}\,dx+C_3\int_{\Omega}|u_\varepsilon|^6\, dx+C_4.
\end{align}
Again recalling \eqref{c3o}, and multiplying \eqref{c3o} by $\frac{16}{5}\chi^2C_p$ on both sides, we see that for any small $\eta, \delta\in (0,1)$,
\begin{align*}
&\frac{8}{5}\chi^2C_p\frac{d}{dt}\int_\Omega |\Delta c_\varepsilon|^{2} dx+\frac{8}{5}\chi^2C_p \int_\Omega |\Delta c_\varepsilon|^{2} dx
+ \frac{6}{5}\chi^2C_p\int_{\Omega}|\nabla\Delta c_\varepsilon|^{2} dx+2\chi^2C_p\int_{\Omega} n_\varepsilon|\Delta c_\varepsilon|^2\,dx\notag
\\
&\le C_p(16\chi\eta)^2\int_{\Omega} \frac{|\nabla n_{\varepsilon}|^{2(p-2)}|D^2 n_\varepsilon|^2}{\rho} dx+\frac{16}{5}\chi^2C_p \delta\int_{\Omega}\frac{|\nabla n_\varepsilon|^{2p}}{\rho^3}\, dx
+ C_5 \int_{\Omega}(|\nabla u_\varepsilon|^3+|u_\varepsilon |^6)\, dx+C_{\eta,\delta}^1,
\end{align*}
Choosing $\delta, \eta$ appropriately small in the above inequality  such that $\frac{16}{5}\chi^2C_p \delta<\frac{\alpha_p}{64}$, $C_p(16\chi\eta)^2<\frac{\alpha_p(p-1)^2}{8p^2}$ gives
\begin{align}
&\frac{8}{5}\chi^2C_p\frac{d}{dt}\int_\Omega |\Delta c_\varepsilon|^{2} dx+\frac{8}{5}\chi^2C_p \int_\Omega |\Delta c_\varepsilon|^{2} dx
+ \frac{6}{5}\chi^2C_p\int_{\Omega}|\nabla\Delta c_\varepsilon|^{2} dx+2\chi^2C_p\int_{\Omega} n_\varepsilon|\Delta c_\varepsilon|^2\,dx\notag
\\
\label{4-36}
&\le \frac{\alpha_p(p-1)^2}{8p^2}\int_{\Omega} \frac{|\nabla n_{\varepsilon}|^{2(p-2)}|D^2 n_\varepsilon|^2}{\rho} dx+\frac{\alpha_p}{64}\int_{\Omega}\frac{|\nabla n_\varepsilon|^{2p}}{\rho^3}\, dx
+ C_5 \int_{\Omega}(|\nabla u_\varepsilon|^3+|u_\varepsilon |^6)\, dx+C_6.
\end{align}
Adding \eqref{4-34}--\eqref{4-36} gives
\begin{align}
&H'(t)+H(t)+\frac{\alpha_p(p-1)^2}{8p^2}\int_\Omega \frac{\big(|\nabla n_\varepsilon|^2+\varepsilon\big)^{p-2}|D^2n_\varepsilon|^2}{\rho}dx +\frac{\alpha_p}{16}\int_\Omega\frac{\big(|\nabla n_\varepsilon|^2+\varepsilon\big)^{p}}{\rho^3}\,dx\notag
\\
&+\int_{\Omega}|\nabla c_\varepsilon|^{2}| D^2 c_\varepsilon|^2 \,dx
+\int_{\Omega} n_\varepsilon|\nabla c_\varepsilon|^{4}\,dx
+ \frac{6}{5}\chi^2C_p\int_{\Omega}|\nabla\Delta c_\varepsilon|^{2} dx+\chi^2C_p\int_{\Omega} n_\varepsilon|\Delta c_\varepsilon|^2\,dx\notag
\\
\label{4-37}
\le & C_7\int_{\Omega}\left(|u_\varepsilon|^{\frac{8p-9}{2p-3}}+|\nabla u_\varepsilon|^3+|u_\varepsilon |^6\right)\, dx+C_8,
\end{align}
where
$$
H(t)=\int_{\Omega}\left(\frac{\big(|\nabla n_\varepsilon|^2+\varepsilon\big)^{\frac{p}{2}}}{\rho}+|\nabla c_\varepsilon|^{4}+\frac{8}{5}\chi^2C_p\int_\Omega |\Delta c_\varepsilon|^{2}\right)dx.
$$
It implies that
\begin{align}
&\sup_{0<t<T_{\max}^\varepsilon} H(t)+\frac{\alpha_p(p-1)^2}{8p^2}\sup_{\tau<t<T_{\max}^\varepsilon}\int_{t-\tau}^t\int_\Omega \frac{\big(|\nabla n_\varepsilon|^2+\varepsilon\big)^{p-2}|D^2n_\varepsilon|^2}{\rho}\,dxds\notag \\
&+\frac{\alpha_p}{16}\sup_{\tau<t<T_{\max}^\varepsilon}\int_{t-\tau}^t\int_\Omega\frac{\big(|\nabla n_\varepsilon|^2+\varepsilon\big)^{p}}{\rho^3}\,dxds\notag
\\
&+\sup_{\tau<t<T_{\max}^\varepsilon}\int_{t-\tau}^t\int_{\Omega}\left(|\nabla c_\varepsilon|^{2}| D^2 c_\varepsilon|^2 +n_\varepsilon|\nabla c_\varepsilon|^{4}
+ \frac{6}{5}\chi^2C_p |\nabla\Delta c_\varepsilon|^{2} +\chi^2C_p n_\varepsilon|\Delta c_\varepsilon|^2\right)\,dxds \notag
\\
\label{4-38}
\le & C_9\sup_{\tau<t<T_{\max}^\varepsilon}\int_{t-\tau}^t\int_{\Omega}\left(|u_\varepsilon|^{\frac{8p-9}{2p-3}}+|\nabla u_\varepsilon|^3+|u_\varepsilon |^6\right)\, dxds+C_{10}.
\end{align}
Recalling \eqref{4-22}, and noting that $\frac{8p-9}{2p-3}<8p-9$, $6<8p-9$, we obtain
$$
C_9\sup_{\tau<t<T_{\max}^\varepsilon}\int_{t-\tau}^t\int_{\Omega}\left(|u_\varepsilon|^{\frac{8p-9}{2p-3}}+|\nabla u_\varepsilon|^3+|u_\varepsilon |^6\right)\, dxds+C_{10}
\le \frac{\alpha_p}{32}\sup_{\tau<t<T_{\max}^\varepsilon}\int_{t-\tau}^t\int_\Omega\frac{\big(|\nabla n_\varepsilon|^2+\varepsilon\big)^{p}}{\rho^3}\,dxds+C_{11}.
$$
Combined with \eqref{4-38}, this completes the proof of Proposition \ref{prop:key}. \hfill $\Box$

\medskip

{\it\bfseries Proof of Theorem \ref{thm:main}.}
Let $p>2$ be fixed. For $p>p^*$ the global boundedness, the first-order bounds in \eqref{1-3} and the $L^2$ estimate for $\partial_tn_\varepsilon$ are contained in \cite{Jin2021JDE}; the genuinely new piece of information is the second-order estimate for the density in \eqref{1-4}. As explained in Remark~\ref{rem:two}, this estimate is obtained by two different arguments according to the size of $p$. We first carry out the argument that applies to $2<p\le3$ (Step~1), then the one that applies to $p>p^*$ (Step~2), and finally pass to the limit $\varepsilon\to0$ (Step~3). Since $(2,3]\cup(p^*,\infty)=(2,\infty)$, the two arguments together cover every $p>2$.

\textbf{Step 1: the estimates for $2<p\le3$.}
The starting point is a uniform $L^3$ bound for the density, which is an immediate consequence of Proposition~\ref{prop:key}. Noticing that
$$
\|\nabla\rho^{\frac{p-1}{p}}(\cdot, t)\|_{L^p}^p=\left(\frac{p-1}{p}\right)^p\int_{\Omega}\frac{|\nabla\rho|^p}{\rho}dx=\left(\frac{p-1}{p}\right)^p\int_{\Omega}\frac{|\nabla n_\varepsilon|^p}{\rho}dx\le \left(\frac{p-1}{p}\right)^p\int_{\Omega}\frac{\big(|\nabla n_\varepsilon|^2+\varepsilon\big)^{\frac{p}{2}}}{\rho}dx,
$$
and that the last integral is bounded by a constant independent of $\varepsilon$ and $t$ by Proposition~\ref{prop:key}, we apply the Sobolev embedding in the form $\|g\|_{L^q}\le C\big(\|\nabla g\|_{L^p}+\|g\|_{L^1}\big)$ to $f:=\rho^{\frac{p-1}{p}}$, and use $\|f\|_{L^1}\le C\|f\|_{L^{\frac{p}{p-1}}}$ together with $\|f\|_{L^{\frac{p}{p-1}}}^{\frac{p}{p-1}}=\int_\Omega\rho\,dx$. This yields
$$
\|\rho(\cdot, t)\|_{L^{\frac{(p-1)q}{p}}}=\|f\|_{L^q}^{\frac{p}{p-1}}\le C\left(\|\nabla f\|_{L^p}+\|f\|_{L^{\frac{p}{p-1}}}\right)^{\frac{p}{p-1}}\le C\left(\|\nabla f\|_{L^p}^{\frac{p}{p-1}}+\|\rho\|_{L^1}\right)\le \tilde C
$$
for any exponent $q$ with
$$
\left\{
\begin{aligned}
&q\le  \frac{3p}{(3-p)}, && \text{for $p<3$},
\\
&q<\infty, &&\text{for $p=3$},
\end{aligned}
\right.
$$
where the last step also uses $\int_\Omega\rho(\cdot,t)\,dx\le \int_\Omega n_0\,dx+|\Omega|$. Since $\frac{p-1}{p}\cdot\frac{3p}{(3-p)}=\frac{3(p-1)}{3-p}>3$ when $2<p<3$, the density is bounded in $L^3(\Omega)$, uniformly in $\varepsilon$ and in $t$.

Combining Lemmas \ref{lem:jin3.4} and \ref{lem:localex} (with $\beta=1$) we therefore infer that
$$
T_{\max}^{\varepsilon}=\infty,
$$
and that there exists a positive constant $C$, independent of $\varepsilon$, such that
$$
\sup_{t>0}\|n_\varepsilon(\cdot, t)\|_{L^\infty}<C.
$$
In particular $\tau=1$ in Proposition~\ref{prop:key}, so that the weighted second-order bounds of Proposition~\ref{prop:key} are global in time. Since now $\rho\le C$, these weighted bounds imply unweighted ones. Precisely, using the pointwise inequality
$|\Delta_p^{\varepsilon}n_\varepsilon|^2\le C_p\big(|\nabla n_\varepsilon|^2+\varepsilon\big)^{p-2}|D^2n_\varepsilon|^2$,
which follows by differentiating $\Delta_p^{\varepsilon}n_\varepsilon=\operatorname{div}\big((|\nabla n_\varepsilon|^2+\varepsilon)^{\frac{p-2}{2}}\nabla n_\varepsilon\big)$ and by $|\nabla n_\varepsilon|\le\big(|\nabla n_\varepsilon|^2+\varepsilon\big)^{\frac12}$, we obtain
\begin{align*}
&\sup_{t>0}\int_t^{t+1}\int_{\Omega}|\Delta_p^{\varepsilon}n_\varepsilon|^2 dxds\le \sup_{t>0}\|\rho(\cdot, t)\|_{L^\infty}\int_t^{t+1}\int_{\Omega}\frac{|\Delta_p^{\varepsilon}n_\varepsilon|^2}{\rho} dxds
\\
&\le C \sup_{t>0}\|\rho(\cdot, t)\|_{L^\infty}\int_t^{t+1}\int_{\Omega}\frac{\left(|\nabla n_\varepsilon|^2+\varepsilon\right)^{p-2}|D^2n_\varepsilon|^2}{\rho} dxds\le \tilde C.
\end{align*}
Similarly, again by the boundedness of $\rho$ and Proposition~\ref{prop:key},
$$
\sup_{t>0}\int_t^{t+1}\int_{\Omega}\left(|\nabla n_\varepsilon|^{2p}+\left(|\nabla n_\varepsilon|^2+\varepsilon\right)^{p-2}|D^2n_\varepsilon|^2\right) dxds\le  C,
$$
since $\rho\le C$ yields
$|\nabla n_\varepsilon|^{2p}\le \big(|\nabla n_\varepsilon|^2+\varepsilon\big)^p\le C\,\frac{\big(|\nabla n_\varepsilon|^2+\varepsilon\big)^p}{\rho^3}$
and
$\big(|\nabla n_\varepsilon|^2+\varepsilon\big)^{p-2}|D^2n_\varepsilon|^2\le C\,\frac{\big(|\nabla n_\varepsilon|^2+\varepsilon\big)^{p-2}|D^2n_\varepsilon|^2}{\rho}$.

With the boundedness of $n_\varepsilon$ at hand, the $L^r$-theory of linear parabolic equations applied to $c_\varepsilon$ and that of the Stokes system applied to $(u_\varepsilon,\pi_\varepsilon)$ give for any $r>1$
$$
\|c_\varepsilon\|_{W_r^{2,1}(Q_1(t))}+\|u_\varepsilon\|_{W_r^{2,1}(Q_1(t))}+\|\nabla \pi_\varepsilon\|_{L^r(Q_1(t))}\le C_r.
$$
Hence $\|u_\varepsilon(\cdot,t)\|_{W^{1,\infty}}\le C$ by the parabolic embedding $W_r^{2,1}(Q_1(t))\hookrightarrow C^{1+\alpha, \frac{\alpha}2}(\overline{Q_1(t)})$, valid for $r>5$; a bound for $\|c_\varepsilon(\cdot,t)\|_{W^{1,\infty}}$ follows from the standard $L^q$-$L^\infty$ estimates for the Neumann heat semigroup, and a bound for $\|A^\beta u_\varepsilon(\cdot,t)\|_{L^2}$ from the corresponding estimates for the Stokes semigroup (see Lemma 3.8 of \cite{Jin2021JDE} for the details). It remains to bound $\partial_tn_\varepsilon$, and for this we simply use the first equation in \eqref{eq:approx}:
$$
\partial_t n_\varepsilon=\Delta_p^{\varepsilon}n_\varepsilon-u_\varepsilon\cdot\nabla n_\varepsilon-\chi\nabla\cdot(n_\varepsilon\nabla c_\varepsilon).
$$
Every term on the right-hand side is controlled on $Q_1(t)$ by the bounds already established: $\Delta_p^{\varepsilon}n_\varepsilon$ and $\Delta c_\varepsilon$ belong to $L^2$, $\nabla n_\varepsilon$ belongs to $L^{2p}$, while $n_\varepsilon$, $u_\varepsilon$ and $\nabla c_\varepsilon$ are bounded. Consequently
$$
\sup_{t>0}\Big(\|c_\varepsilon(\cdot, t)\|_{W^{1,\infty}}+\|u_\varepsilon(\cdot, t)\|_{W^{1,\infty}}+\|A^\beta u_\varepsilon(\cdot, t)\|_{L^2}\Big)+\sup_{t>0}\|\partial_t n_\varepsilon\|_{L^2(Q_1(t))}\le C,
$$
which, together with the bounds obtained above, is precisely \eqref{1-3}--\eqref{1-4} for $2<p\le3$, uniformly in $\varepsilon$.

\textbf{Step 2: the second-order estimate for $p>p^*$.}
In this range no new first-order information is required, and the second-order estimate for the density is obtained without any density weight. The global boundedness, the first-order bounds in \eqref{1-3}, the $L^r$-bounds for $c_\varepsilon, u_\varepsilon, \nabla\pi_\varepsilon$ (which follow from the boundedness of $n_\varepsilon$ exactly as in Step~1) and the $L^2$ estimate for $\partial_tn_\varepsilon$ are all available from \cite{Jin2021JDE}; in particular, Lemma 3.10 there gives
$$
\sup_{t>0}\int_t^{t+1}\int_\Omega|\partial_t n_\varepsilon|^2 dxds\le C.
$$
To extract the second-order estimate for the density, we return to the first equation in \eqref{eq:approx}:
$$
\Delta_p^{\varepsilon} n_\varepsilon=\partial_t n_\varepsilon+u_\varepsilon\cdot\nabla n_\varepsilon
+\chi\nabla\cdot(n_\varepsilon\nabla c_\varepsilon).
$$
Combining the estimate for $\partial_tn_\varepsilon$ with the bounds recalled above for $n_\varepsilon$, $\nabla n_\varepsilon$, $c_\varepsilon$ and $u_\varepsilon$ (recall that $c_\varepsilon\in W^{2,1}_r(Q_1(t))$ for every $r>1$, so that $\Delta c_\varepsilon\in L^2$), we obtain
$$
\sup_{t>0}\|\Delta_p^{\varepsilon} n_\varepsilon\|_{L^2(Q_1(t))}\leq C.
$$
To transfer this information from $\Delta_p^{\varepsilon}n_\varepsilon$ to the Hessian, we employ an unweighted Rellich identity. More precisely, with $A_\varepsilon$ and $Z_\varepsilon$ defined as in \eqref{eq:MZ}, the same computation as in the proof of Lemma \ref{lem:rellich} gives, for any $p>2$,
$$
\int_\Omega(\Delta_p^\varepsilon n_\varepsilon)^2\,dx
=\int_\Omega |Z_\varepsilon|^2\,dx+\int_{\partial\Omega}\mathrm{II}(A_\varepsilon,A_\varepsilon)\,ds.
$$
The boundary contribution is now controlled by adapting the trace inequality and the Sobolev interpolation inequality; no weight $\rho$ occurs in these arguments, which is what makes this route simpler than the one of Step~1. Consequently,
$$
\sup_{t>0}
\int_t^{t+1}\int_\Omega |Z_\varepsilon|^2 \,dx\,ds\leq C,
$$
and hence, by the spectral bounds for $M_\varepsilon$ (Lemmas \ref{lem:Mv} and \ref{lem:eig}),
$$
\sup_{t>0}\int_t^{t+1}\int_\Omega \big(|\nabla n_\varepsilon|^2+\varepsilon\big)^{p-2}|D^2n_\varepsilon|^2 \,dx\,ds\leq C.
$$
The above estimates also give a uniform bound for $\nabla n_\varepsilon$ in $L^{2p}$. Indeed, applying \eqref{ineq-2} with $r=2p-2$ to $n_\varepsilon(\cdot,t)$ and integrating in $t$, we obtain
$$
\sup_{t>0}\int_t^{t+1}\int_\Omega |\nabla n_\varepsilon|^{2p}\,dxds\le (2p-2+\sqrt 3)^2\sup_{t>0}\|n_\varepsilon(\cdot,t)\|_{L^\infty}^2\sup_{t>0}\int_t^{t+1}\int_{\Omega}|\nabla n_\varepsilon|^{2(p-2)}|D^2n_\varepsilon|^2\,dxds\le \tilde C.
$$
This establishes \eqref{1-3} and \eqref{1-4} for $p>p^*$; together with Step~1, we conclude that both hold throughout the whole range $p>2$.

\medskip

\textbf{Step 3: passage to the limit.}
By Steps~1 and~2, for every fixed $p>2$ the bounds in \eqref{1-3} and \eqref{1-4} hold uniformly in $\varepsilon$ and in $t$. With these uniform bounds at hand, the Aubin--Lions lemma, applied to the equations in \eqref{eq:approx} on $Q_T=\Omega\times(0,T)$, allows us to pass to the limit $\varepsilon\to 0$ along a subsequence, which we do not relabel, so that for any given $Q_T=\Omega\times(0, T)$,
\begin{align*}
&u_\varepsilon, c_\varepsilon \to  u, c, && \text{ uniformly in $Q_T$};
\\
&\nabla c_\varepsilon \to \nabla c, && \text{in $L^r(Q_T)$ for any $r>1$};
\\
&n_\varepsilon \to n, && \text{in $L^q(Q_T)$ for any $q>1$};
\\
&n_\varepsilon \xrightarrow{w^*} n, && \text{in $L^\infty(Q_T)$};
\\
&\nabla n_\varepsilon \to \nabla n, && \text{in $L^q(Q_T)$ for any $1<q<2p$};
\\
&\partial_t n_\varepsilon  \xrightarrow{w}\partial_t n, \,\Delta_p^{\varepsilon} n_\varepsilon \xrightarrow{w} \Delta_p n, && \text{weakly in $L^2(Q_T)$};
\\
&\partial_t c_\varepsilon \xrightarrow{w} \partial_t c, \, \partial_t u_\varepsilon \xrightarrow{w}\partial_t u, \,\nabla\pi_\varepsilon\xrightarrow{w}\nabla\pi, \,\Delta u_\varepsilon \xrightarrow{w} \Delta u, && \text{in $L^r(Q_T)$ for any $r>1$};
\\
&\nabla\Delta c_\varepsilon \xrightarrow{w} \nabla\Delta c && \text{in $L^2(Q_T)$ };
\\
&\Delta c_\varepsilon \to \Delta c, && \text{in $L^r(Q_T)$ for any $r>1$};
\\
& \nabla\cdot(n_\varepsilon\nabla c_\varepsilon)\to \nabla\cdot(n\nabla c), && \text{in $L^2(Q_T)$}.
\end{align*}
With these convergences, every term in \eqref{eq:approx} passes to the limit in the sense of distributions, so that $(n,c,u,\pi)$ solves \eqref{eq:main}--\eqref{eq:ic}. Moreover $(n,c,u,\pi)\in\mathcal X_1\times\mathcal X_2\times\mathcal X_3\times \mathcal X_4$, and \eqref{1-3}--\eqref{1-4} follow by weak lower semicontinuity. This completes the proof of Theorem \ref{thm:main}.
 \hfill $\Box$

\end{document}